%% file: 1quadruples.tex
\documentclass[a4paper,12pt,reqno]{amsart}

\usepackage{graphicx}
\usepackage{amsmath}
\usepackage{amssymb}
\usepackage{amsfonts}
\usepackage{epsfig}

\makeatletter

\newtheorem{theorem}{Theorem}[section]

 \@addtoreset{equation}{section}
 \@addtoreset{figure}{section}
 \@addtoreset{table}{section}
 \@addtoreset{theorem}{subsection}

\newtheorem{corollary}[theorem]{Corollary}
\newtheorem{remark}[theorem]{Remark}
\newtheorem{proposition}[theorem]{Proposition}
\newtheorem{lemma}[theorem]{Lemma}

\newtheorem{conjecture}[theorem]{Conjecture}
\def\PerfProof{{\it Proof.\ }}

\def\rep{{\rm rep}}
\def\longarr#1#2{{\buildrel{#1} \over {\hbox to
#2pt{\rightarrowfill}}}}
\def\indownarrow{{\raise5.5pt\hbox{$\scriptstyle\cap$}}\!\!\!\hspace{-0.15em}\downarrow}
\def\inuparrow{{\raise-3.5pt\hbox{$\scriptstyle\cup$}}\!\!\!\hspace{-0.15em}\uparrow}
\def\longdownarrow{{\raise5.5pt\hbox{$\scriptstyle{\arrowvert}$}}\!\!\!\hspace{-0.15em}\downarrow}
\def\longuparrowr{{\raise5.5pt\hbox{$\scriptstyle{\arrowvert}$}}\!\!\!\hspace{-0.15em}\uparrow}
\def\longuparrow{{\raise-3.5pt\hbox{$\scriptstyle{\arrowvert}$}}\!\!\!\hspace{-0.15em}\uparrow}

\makeatother

\makeindex

\begin{document}

\title[Quadruples, admissible elements and Herrmann's endomorphisms]
{Quadruples, admissible elements and Herrmann's endomorphisms}
         \author{Rafael Stekolshchik}
         \thanks{{\large email: rs2@biu.013.net.il}}

\date{}

\begin{abstract}
  The notion of a {\it perfect element} of
 the free modular lattice $D^r$ generated by $r \geq 1$ elements
 is introduced by Gelfand and Ponomarev in \cite{GP74}.
  The classification of such elements in the lattice $D^4$
 (and also in $D^r$, where $r > 4$)
 is given in \cite{GP74}, \cite{GP76}, \cite{GP77}.

 {\it Admissible elements} are the building blocks used
 for the construction of perfect elements.
 In \cite{GP74}, Gelfand and Ponomarev construct
 admissible elements of the lattice $D^r$ recurrently.
 We suggest a direct method for creating
 admissible elements in $D^4$ and show that these elements coincide
 modulo linear equivalence with admissible elements constructed
 by Gelfand and Ponomarev. Admissible sequences and admissible elements
 for $D^4$ form $8$ classes and possess some periodicity.

 The main aim of the present paper is to find a connection between
 admissible elements in $D^4$ and {\it Herrmann's endomorphisms}
 $\gamma_{1i}$. Herrmann in \cite{H82} constructed perfect elements $s_n$, $t_n$,
 $p_{i,n}$ in $D^4$ by means of some {\it endomorphisms}
 $\gamma_{ij}$ and showed that these perfect elements coincide with
 the Gelfand-Ponomarev perfect elements modulo linear equivalence.
 We show that the admissible elements in $D^4$ are also obtained by
 means of Herrmann's endomorphisms $\gamma_{ij}$. Herrmann's
 endomorphism $\gamma_{ij}$ and the {\it elementary map} of
 Gelfand-Ponomarev $\phi_i$ act, in a sense, in opposite
 directions, namely the endomorphism $\gamma_{ij}$ adds the index
 to the start of the admissible sequence, and the elementary map
 $\phi_i$ adds the index to the end of the admissible sequence, see
 (\ref{map_phi_adds}) and (\ref{map_Herrmann_adds}).
\end{abstract}



\maketitle

 \tableofcontents
 \listoffigures
 \listoftables
\newpage

\input 2outline_adm.tex
\input 3admiss_D4.tex
\input 4endom_D4.tex

\renewcommand{\appendixname}{}


\end{document}

%% file: 2outline_adm.tex
\section{\sc\bf Introduction}
We denote by $D^4$ the free modular lattice with four generators.
Indecomposable representations of $D^4$ were classified by
Nazarova in \cite{Naz67} and Gelfand-Ponomarev in
\cite{GP70}\footnote{Equivalently, instead of representations of
$D^4$, one can speak about {\it quadruple of subspaces},
\cite{GP70}.}.

In \cite{GP74} Gelfand-Ponomarev constructed {\it perfect
elements} in $D^4$. We recall that an element $a \in D$ of a
modular lattice $D$ is perfect, if for each finite dimension
indecomposable $K$-linear representation $\rho_X : L \rightarrow
\mathfrak{L}(X)$ over any field $K$, the image $\rho_X(a)
\subseteq X$ of $a$ is either zero, or $\rho_X(a) = X$, where
$\mathfrak{L}(X)$ is the lattice of all vector $K$-subspaces of
$X$.

{\it Admissible elements\footnote{The definition of notion {\it
admissible} will be given shortly.}} are the building blocks in
the construction of perfect elements, see \cite{GP74},
\cite{GP76}, and also \cite{St04}. Making use of certain elegant
endomorphisms in $D^4$, Herrmann in \cite{H82} also constructed
perfect elements. The main purpose of this article is to obtain
{\it a connection between admissible elements and Herrmann's
endomorphisms}.

First, in this section, we outline the idea of admissible elements
and the construction of the distributive sublattice of perfect
polynomials \cite{GP74}. We think today, that apart from being
helpful in construction of perfect elements, the admissible
elements are interesting in themselves.

 Further, we outline properties of
 admissible elements obtained in this work:
 {\it finite classification}, {\it $\varphi-$homomorphism},
 {\it reduction to atomic elements} and {\it periodicity}, see
 Section \ref{adm_steps}. In Section \ref{sect_adm_Herrmanm_endom}
 we outline the connection between admissible elements
 and Herrmann's endomorphisms \cite{H82}.
 Not only perfect elements in $D^4$, also
 the admissible elements $e_\alpha$ and $f_{\alpha0}$ in $D^4$
 can be obtained by Herrmann's endomorphisms, see
 Theorem \ref{theor_adm_Herrmann}.

 In Section \ref{sect_adm_seq_D4}, we construct admissible elements
 in $D^4$, by applying the technique
 developed in \cite{St04} for the modular lattice $D^{2,2,2}$.

 In Section \ref{sect_Herrmann} we study Herrmann's endomorphisms
 used in his construction of the perfect elements in $D^4$.
 We show {\it how the admissible elements introduced in Section \ref{sect_adm_seq_D4}
 are constructed by means of Herrmann's endomorphisms} $\gamma_{ij}$.
 Set of admissible elements $e_\alpha$ (resp. $f_{\alpha0}$)
 (described by Table \ref{table_adm_elem_D4}) coincide with the set
 \begin{equation}
    \begin{split}
       \{ e_i{a}_r^{jl}{a}_s^{lk}{a}_t^{kj} \}, \text{ resp. }
       \{ e_i{a}_r^{jl}{a}_s^{lk}{a}_t^{kj}
           (e_i{a}_{t+1}^{jk}+ {a}_{s+1}^{kl}a_{r-1}^{jl}) \},
   \end{split}
 \end{equation}
where ${a}_r^{jl}$ are {\it atomic elements} (\ref{def_atomic}),
$r,s,t \in \mathbb{Z}_+$ and $\{i,j,k,l\} = \{1,2,3,4\}$, see
Proposition \ref{prop_gen_form_D4}.

 Herrmann's endomorphism $\gamma_{ij}$ and the elementary map of
 Gelfand-Ponomarev $\phi_i$ act, in a sense, in opposite
 directions, namely the elementary map $\phi_i$ adds the index to
 the end of the admissible sequence:
 \begin{equation}
  \label{map_phi_adds}
       \varphi_i\rho_{X^+}(z_\alpha) = \rho_X(z_{i\alpha}),
 \end{equation}
 see Theorem \ref{th_adm_classes_D4}, and the endomorphism
 $\gamma_{ik}$ adds the index to the start of the admissible
 sequence:
 \begin{equation}
  \label{map_Herrmann_adds}
       \gamma_{ik}(z_{{\alpha}{k}}) = z_{\alpha{ki}},
 \end{equation}
see Theorem \ref{relat_rho_ik}.

 Admissible elements $e_\alpha$ (similarly, $f_{\alpha0}$)
     are obtained by means of Herrmann's endomorphisms as follows:
  \begin{equation}
           \gamma^t_{ij}\gamma^s_{ik}\gamma^r_{il}(e_i) =
             e_i{a}^{kl}_t{a}^{jl}_s{a}^{kj}_r,
  \end{equation}
 where $r,s,t \in \mathbb{Z}_+$ and $\{i,j,k,l\} = \{1,2,3,4\}$,
 see Theorem \ref{theor_adm_Herrmann}.

\subsection{The idea of admissible elements of Gelfand and
Ponomarev }
  \label{idea_adm}
  \index{Gelfand-Ponomarev!- elementary maps $\varphi_i$}
  \index{Gelfand-Ponomarev!- admissible elements}

Here, we describe, omitting some details, the idea of constructing
the elementary maps and admissible elements of Gelfand-Ponomarev.
As we will see below, the admissible elements {\it grow}, in a
sense, from generators of the lattice or from the lattice's unity.

Following Bernstein, Gelfand and Ponomarev \cite{BGP73}, given a
modular lattice $L$ and a field $K$, we use the Coxeter functor
$$
    \Phi^+ : \rep_K{L} \longrightarrow \rep_K{L},
$$
see \cite[App. A]{St04} for details. Given a representation
 $\rho_X : L \longrightarrow \mathcal{L}(X)$ in $\rep_K{L}$,
 we denote by
$$
   \rho_{X^+} :
     L \longrightarrow \mathcal{L}(X^+)
$$
the image representation $\Phi^+\rho_X$ in $\rep_K{L}$ under the
functor
$\Phi^+$.  

{\it Construction of the elementary map $\varphi$}.
 Let there exist a map $\varphi$
 mapping every subspace
 $A \in \mathcal{L}(X^+)$ to some subspace
 $B \in \mathcal{L}(X)$:
\begin{equation}
  \label{phi_def_0}
  \begin{array}{c}
    \varphi: \mathcal{L}(X^+) \longrightarrow \mathcal{L}(X), \vspace{2mm}\\
    \varphi{A} = B, \quad \text{ or } \quad
     A\stackrel{\varphi}{\longmapsto}B.
  \end{array}
\end{equation}
 It turns out that for many elements $a \in L$ there
 exists an element $b \in L$ (at least in $L = D^4$ and $D^{2,2,2}$) such that
\begin{equation}
 \label{phi_rho_ab}
   \varphi(\rho_{X^+}(a)) = \rho_X(b)
\end{equation}
for every pair of representations $(\rho_X, \rho_{X^+})$, where
$\rho_{X^+} = \Phi^+\rho_{X}$. The map $\varphi$ in
(\ref{phi_rho_ab}) is constructed in such a way, that $\varphi$ do
not depend on $\rho_X$. Thus, we can write
\begin{equation}
  \label{rel_elem_map_0}
   a \stackrel{\varphi}{\longmapsto} b.
\end{equation}
In particular, eq. (\ref{phi_rho_ab}) and (\ref{rel_elem_map_0})
are true for admissible elements, whose definition will be given
later, see Section 2. Eq. (\ref{phi_rho_ab}) is the main property
characterizing the admissible elements.

For $D^4$, Gelfand and Ponomarev constructed $4$ maps of form
(\ref{phi_def_0}): $\varphi_1, \varphi_2, \varphi_3, \varphi_4$
and called them {\it elementary maps}.  Given a sequence of
relations
 \index{admissible elements! - growing from generators}
\begin{equation}
  \label{rel_elem_map_1}
   a_{i_1} \stackrel{\varphi_{i_1}}{\longmapsto}
   a_{i_2} \stackrel{\varphi_{i_2}}{\longmapsto} a_{i_3}
               \stackrel{\varphi_{i_3}}{\longmapsto}
   \dots
               \stackrel{\varphi_{i_{n-2}}}{\longmapsto}
   a_{i_{n-1}} \stackrel{\varphi_{i_{n-1}}}{\longmapsto} a_n,
\end{equation}
we say that elements
 $a_{i_2},a_{i_3}, \dots, a_{i_{n-2}}, a_{i_{n-1}}, a_{i_n}$
 {\it grow} from the element $a_{i_1}$.
 By abuse of language, we say that elements growing from
 generators $e_i$ or unity $I$ are {\it admissible elements}, and
 the corresponding sequence of indices $\{i_{n}i_{n-1}i_{n-2}\dots{i}_2{i}_1\}$
 is said to be the {\it admissible sequence}.
 For details in the cases $D^{2,2,2}$ and $D^4$,
 see \cite{St04}.

\subsection{Reduction of the admissible elements to atomic elements}
  \label{adm_steps}
Here, we briefly give steps of creation of the admissible elements
in this work.

 \index{key property of admissible sequences}
 \index{admissible sequences! - key property}

\begin{figure}[h]
\centering
\includegraphics{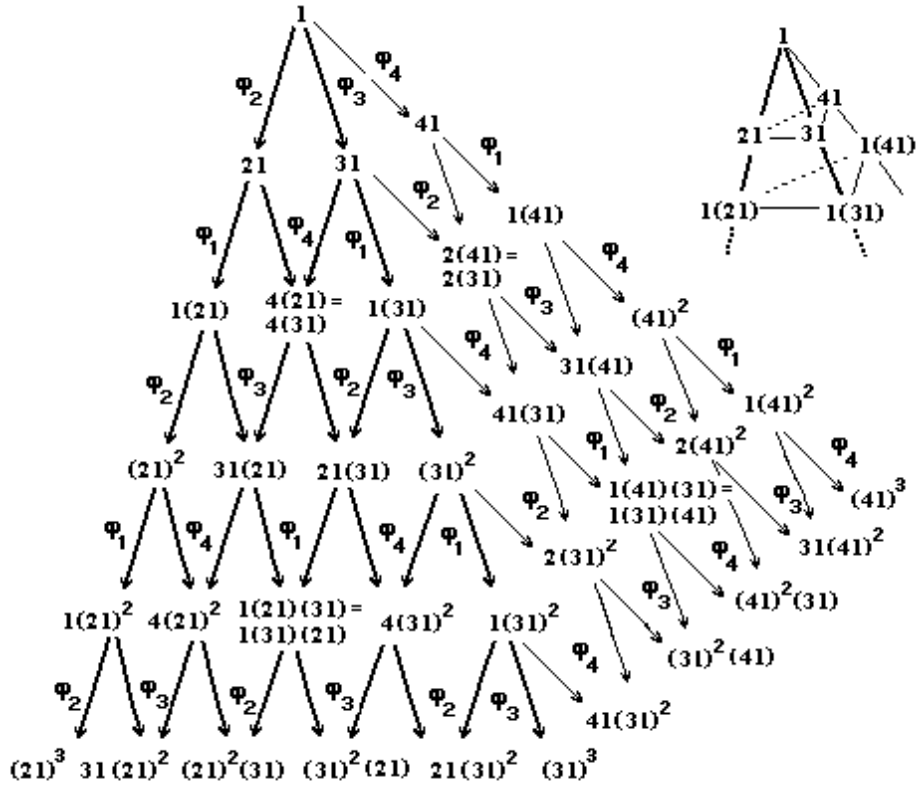}
\caption{\hspace{3mm}Action of maps $\varphi_i$ on the admissible
sequences for $D^4$}
\label{pyramid_D4}
\end{figure}

 {\it A) Finite classification}.
 The admissible elements (and admissible sequences) can be
 reduced to a finite number of classes, see \cite{St04}.
 Admissible sequences for $D^4$ are depicted
 on Fig. \ref{pyramid_D4}.

 The key property of admissible sequences allowing to do this
 classification is the following relation between maps $\varphi_i$.

  \underline{For $D^{2,2,2}$} (see \cite[Prop. 2.10.2]{St04}):
\begin{equation*}
  \begin{split}
    iji & = iki, \\
    \varphi_i\varphi_j\varphi_i & = \varphi_i\varphi_k\varphi_i,
      \text{ where } \{i,j,k\} = \{1,2,3\}.
  \end{split}
\end{equation*}

 \underline{For $D^4$} (see \cite[Prop. 4.6.1]{St04}):
\begin{equation*}
  \begin{split}
    ikj & = ilj, \\
    \varphi_i\varphi_k\varphi_j & = \varphi_i\varphi_l\varphi_j,
    \text{ where } \{i,j,k,l\} = \{1,2,3,4\}.
  \end{split}
\end{equation*}

 {\it B) $\varphi_i-$homomorphism}.
 In a sense, the elementary maps $\varphi_i$
 are homomorphic with respect to admissible elements. More exactly,
 introduce a notion of $\varphi_i-$homomorphic elements. Let
\begin{equation}
  \label{phi_homom_0}
   a \stackrel{\varphi_i}{\longmapsto} \tilde{a}
    \quad \text{ and }\quad
   p \stackrel{\varphi_i}{\longmapsto} \tilde{p}.
\end{equation}
The element $a \in L$ is said to be {\it $\varphi_i-$homomorphic},
if
\begin{equation}
  \label{phi_homom_1}
   ap \stackrel{\varphi_i}{\longmapsto} \tilde{a}\tilde{p}
\end{equation}
for all $p \in L$.

 \index{atomic elements}

 {\it C) Reduction to atomic elements}.
 We permanently apply the following mechanism of creation of admissible elements
from more simple elements. Let $\varphi_i$ be any elementary map,
and
\begin{equation}
  \label{adm_chains}
    a \stackrel{\varphi_i}{\longmapsto} \tilde{a}, \quad
    b \stackrel{\varphi_i}{\longmapsto} \tilde{b}, \quad
    c \stackrel{\varphi_i}{\longmapsto} \tilde{c}, \quad
    p \stackrel{\varphi_i}{\longmapsto} \tilde{p}.
\end{equation}
Suppose $abcp$ is any admissible element and
 $a, b, c$ are $\varphi_i-$homomorphic.
By means of relation (\ref{phi_homom_1}) we construct new
admissible element
\begin{equation*}
  \tilde{a}\tilde{b}\tilde{c}\tilde{p}, \quad
\end{equation*}
The elements $a, b, c$ are called {\it atomic}.  For definition of
atomic elements in $D^{2,2,2}$, see \cite{St04}.

For the case of $D^4$, {\it atomic} lattice polynomials
$a_n^{ij}$, where $i,j \in \{1,2,3,4\}$, $n \in \mathbb{Z}_+$, are
defined as follows
\begin{equation}
 \label{def_atomic}
  a_n^{ij} =
 \begin{cases}
   a_n^{ij} = I \qquad \qquad \qquad \qquad \qquad \hspace{3mm} \text{ for } n = 0,\\
   a_n^{ij} = e_i + e_j{a}_{n-1}^{kl} = e_i + e_j{a}_{n-1}^{lk}
        \hspace{1.7mm} \text{ for } n \geq 1,
 \end{cases}
\end{equation}
where $\{i,j,k,l\}$ is the permutation of $\{1,2,3,4\}$.

\begin{figure}[h]
\includegraphics{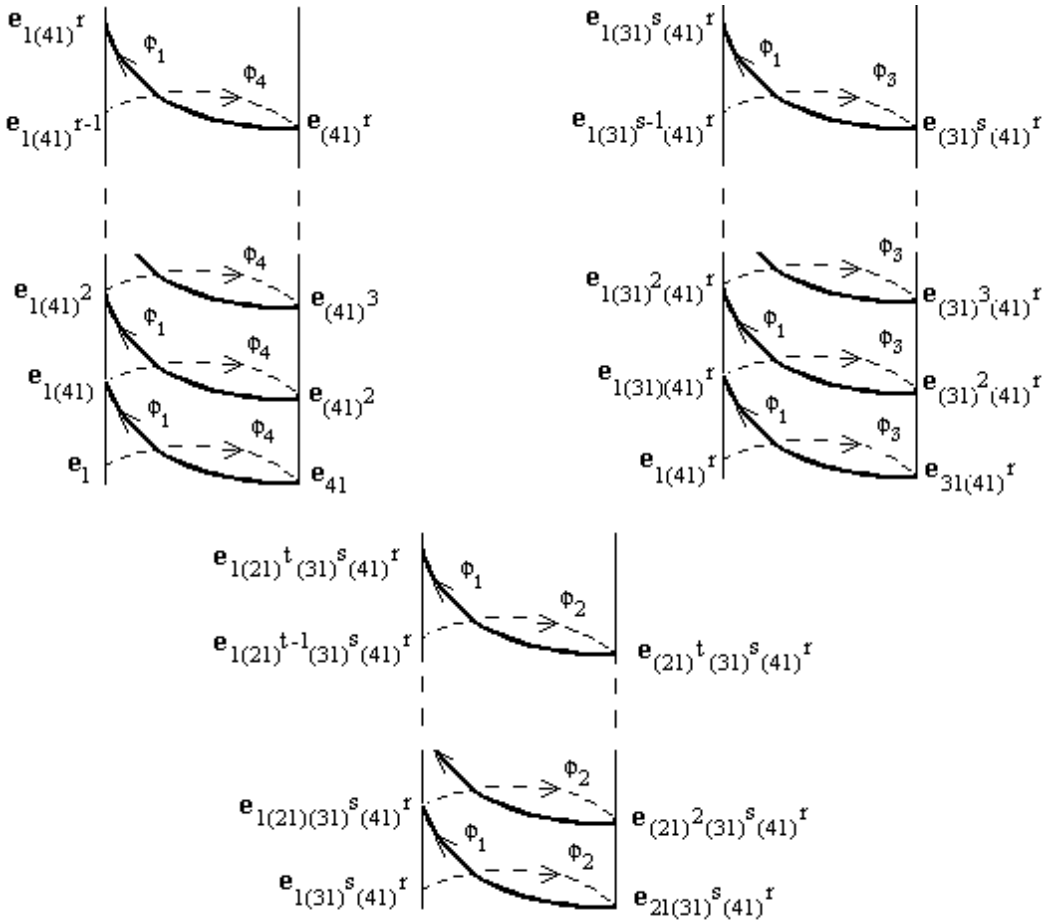}
\caption{\hspace{3mm}Periodicity of admissible elements}
\label{periodicity}
\end{figure}

{\it D) Periodicity}. Admissible sequences obtained by reducing in
heading A) possess some {\it periodicity}. The corresponding
admissible elements are also periodic. On Fig. \ref{periodicity}
we see some examples of admissible elements for $D^4$. Two
vertical lines on each of three cylinders on Fig.
\ref{periodicity} correspond to six series of inclusions of
admissible elements:
\begin{equation}
 \label{period_inclusion}
  \begin{split}
  & \dots \subseteq e_{1(41)^r}
         \subseteq e_{1(41)^{r-1}}
         \subseteq \dots
         \subseteq e_{1(41)}
         \subseteq e_{1}, \vspace{2mm} \\
 & \dots \subseteq e_{(41)^r}
         \subseteq e_{(41)^{r-1}}
         \subseteq \dots
         \subseteq e_{(41)}, \vspace{2mm} \\
 & \dots \subseteq e_{1(31)^s(41)^r}
         \subseteq e_{1(31)^{s-1}(41)^r}
         \subseteq \dots
         \subseteq e_{1(31)(41)^r}
         \subseteq e_{1(41)^r}, \vspace{2mm} \\
 & \dots \subseteq e_{(31)^s(41)^r}
         \subseteq e_{(31)^{s-1}(41)^r}
         \subseteq \dots
         \subseteq e_{(31)(41)^r}
         \subseteq e_{(41)^r}, \\
 & \dots \subseteq e_{1(21)^t(31)^s(41)^r}
         \subseteq e_{1(21)^{t-1}(31)^s(41)^r}
         \subseteq \dots
         \subseteq e_{1(21)(31)^s(41)^r}
         \subseteq e_{1(31)^s(41)^r}, \vspace{2mm} \\
 & \dots \subseteq e_{(21)^t(31)^s(41)^r}
         \subseteq e_{(21)^{t-1}(31)^s(41)^r}
         \subseteq \dots
         \subseteq e_{(21)(31)^s(41)^r}
         \subseteq e_{(31)^s(41)^r}. \vspace{2mm} \\
 \end{split}
\end{equation}

  Three cylindric helices on Fig. \ref{periodicity} correspond
to relations
\begin{equation}
  \label{helix_adm}
  \begin{array}{lll}
  \vspace{0.3mm} \\
   e_{1(41)^{r-1}} \stackrel{\varphi_4}{\longmapsto} e_{(41)^r},
    &  & e_{1(31)^{s-1}(41)^r} \stackrel{\varphi_3}{\longmapsto} e_{(31)^s(41)^r}, \\
   e_{(41)^r} \stackrel{\varphi_1}{\longmapsto} e_{1^(41)^r},
    &  & e_{(31)^s(41)^r}  \stackrel{\varphi_1}{\longmapsto}  e_{1(31)^s(41)^r}, \\
    & e_{1(21)^{t-1}(31)^s(41)^r}
      \stackrel{\varphi_2}{\longmapsto} e_{(21)^t(31)^s(41)^r}, \\
    & e_{(21)^t(31)^s(41)^r}  \stackrel{\varphi_1}{\longmapsto}
      e_{1(21)^t(31)^s(41)^r}. \\
      \vspace{0.3mm} \\
  \end{array}
\end{equation}
Relation (\ref{helix_adm}) is a particular case of Theorem
\ref{th_adm_classes_D4} (theorem on admissible elements;
 compare also with the theorem on admissible elements for $D^{2,2,2}$
 \cite[Th. 2.12.1]{St04}).

\subsection{Direct construction of admissible elements}

Admissible elements for $D^r$ in \cite{GP74}, \cite{GP76} are
built recurrently in the length of multi-indices named {\it
admissible sequences}. In this work we suggest a direct method for
creating admissible elements. For $D^{2,2,2}$ (resp. $D^4$), the
admissible sequences and admissible elements form 14 classes
(resp. 8 classes) and possess some periodicity properties, see
Section 1.7.1, Tables 2.2., 2.3, (resp. Section 1.7.2, Tables 4.1,
4.3) from \cite{St04}.

For the definition of admissible polynomials in $D^4$, see
 Table \ref{table_adm_elem_D4}.

For the definition of admissible polynomials in $D^4$ due to
Gelfand and Ponomarev and examples obtained from this definition,
see Section \ref{coinc_GP_D4}.

\subsection{Eight types of admissible sequences}
 \label{adm_D4}

The admissible sequence for $D^4$ is defined as follows. Consider
a finite sequence of indices $s = i_n\dots{i}_1$, where $i_p \in
\{1,2,3\}$. The sequence $s$ is said to be {\it admissible} if
\begin{equation*}
 \begin{split}
  & \text{ (a) Adjacent indices are distinct ($i_p \neq i_{p+1}$)}. \\
  & \text{ (b) In each subsequence $ijl$, we can replace index
    $j$ by $k$. In other words: } \\
  & \dots{ijl}\dots = \dots{ikl}\dots , \text{ where all indices } i,j,k,l
 \text{ are distinct}.
 \end{split}
\end{equation*}

Any admissible sequence with $i=1$ for $D^4$ may be transformed to
one of the next $8$ types (Proposition \ref{full_adm_seq_D4},
Table \ref{table_admissible_ExtD4}):
\begin{equation}
\begin{split}
 & 1) ~(21)^t(41)^r(31)^s = (21)^t(31)^s(41)^r, \\
 & 2) ~(31)^t(41)^r(21)^s = (31)^t(21)^s(41)^r, \\
 & 3) ~(41)^t(31)^r(21)^s = (41)^t(21)^s(31)^r, \\
 & 4) ~1(41)^t(31)^r(21)^s = 1(31)^t(41)^s(21)^r =
 1(21)^t(31)^s(41)^r,\\
 & 5) ~2(41)^r(31)^s(21)^t = 2(31)^{s+1}(41)^{r-1}(21)^t, \\
 & 6) ~3(41)^r(21)^s(31)^t = 3(21)^{s+1}(41)^{r-1}(31)^t, \\
 & 7) ~4(21)^r(31)^s(41)^t = 3(31)^{s+1}(21)^{r-1}(41)^t, \\
 & 8) ~(14)^r(31)^s(21)^t = (14)^r(21)^{t+1}(31)^{s-1} =
       (13)^s(41)^r(21)^t = \\
 &     ~(13)^s(21)^{t+1}(41)^{r-1} = (12)^{t+1}(41)^r(31)^{s-1} =
       (12)^{t+1}(31)^{s-1}(41)^r. \\
 \end{split}
 \end{equation}

\subsection{Admissible elements in $D^4$ and Herrmann's endomorphisms}
  \label{sect_adm_Herrmanm_endom}
 C.~Herrmann introduced in \cite{H82} the commuting endomorphisms $\gamma_{ik}$
  ($i,k = 1,2,3$) in $D^4$ and used them in the construction of the perfect
  elements, see Section \ref{sect_Herrmann}.

It turned out, that endomorphisms $\gamma_{ik}$ are also closely
connected with admissible elements. First of all, the endomorphism
$\gamma_{ik}$ acts on
 the admissible element $e_{\alpha{k}}$ such that
 \begin{equation*}
    \gamma_{ik}(e_{\alpha{k}}) = e_{\alpha{ki}},
 \end{equation*}
  There is a {\it similarity between the action of Herrmann's
 endomorphism $\gamma_{ik}$ and the action of the elementary map of Gelfand-Ponomarev $\phi_i$}.
 The endomorphism $\gamma_{ik}$ and the elementary map $\phi_i$ act,
 in a sense, in opposite directions,
 namely the endomorphism $\gamma_{ik}$ adds the index
 to the \underline{start} of the admissible sequence
 (see Theorem \ref{relat_rho_ik}), and the elementary map $\phi_i$ adds
 the index to the \underline{end} (see Theorem \ref{th_adm_classes_D4}).

 Further, every admissible element $e_\alpha$ (resp. $f_{\alpha0}$)
 in $D^4$ has an unified form
 $e_1{a}^{34}_t{a}^{24}_s{a}^{32}_r$, see Table
 \ref{tab_gen_form_D4_e}, and
 every admissible element $e_\alpha$ (resp. $f_{\alpha0}$)
 is obtained by means of the sequence of Herrmann's
 endomorphisms as follows:
 \begin{equation*}
      \gamma^t_{12}\gamma^s_{13}\gamma^r_{14}(e_1) =
             e_1{a}^{34}_t{a}^{24}_s{a}^{32}_r,
 \end{equation*}
 and, respectively,
\begin{equation*}
      \gamma^t_{12}\gamma^s_{13}\gamma^r_{14}(f_{10}) =
             e_1{a}^{34}_t{a}^{24}_s{a}^{32}_r
               (e_1a^{34}_{t+1} + a^{24}_{s+1}a^{32}_{r-1}),
 \end{equation*}
 see Theorem \ref{theor_adm_Herrmann}.

\subsection{Examples of admissible elements}
  \label{examples_D4}
Admissible sequences and admissible elements are taken from Table
 \ref{table_adm_elem_D4}.

  \underline{For $n=1$}: \indent $\alpha = 1$,
  (Table \ref{table_adm_elem_D4}, Line $G11$, $r = 0, s = 0, t = 0$). We have
\begin{equation*}
 \begin{split}
 & e_1 = e_1 \quad \text{(admissible element $e_1$ coincides with
  generator $e_1$)}, \\
 & f_{10} = e_1(e_2 + e_3 + e_4) \subseteq e_1.
 \end{split}
\end{equation*}

  \underline{For $n=2$}: \indent $\alpha = 21$,
  (Table \ref{table_adm_elem_D4}, Line $G21$, $r = 0,s = 0, t = 1$). We have
\begin{equation*}
 \begin{split}
  e_{21} = &  e_2(e_3 + e_4), \\
  f_{210} = & e_2(e_3 + e_4)(e_4 + e_3(e_1 + e_2) + e_1) = \\
 &         e_2(e_3 + e_4)(e_4 + e_2(e_1 + e_3) + e_1) =  \\
 &         e_2(e_3 + e_4)(e_2(e_1 + e_4) + e_2(e_1 + e_3)) \subseteq e_{21}.
 \end{split}
\end{equation*}

  \underline{For $n=3$}: \indent Consider two admissible sequences:
  $\alpha = 121$ and  $\alpha = 321 = 341$. \vspace{2mm}

  1) $\alpha = 121$,
  (Table \ref{table_adm_elem_D4}, Line $G11$, $r = 0,s = 0, t = 1$). We have
\begin{equation*}
 \begin{split}
   e_{121} = & e_1{a}^{34}_2 = e_1(e_3 + e_4(e_1 + e_2)) =  \\
   & e_1(e_3(e_1 + e_2) + e_4(e_1 + e_2)), \\
  f_{1210} = & e_{121}(e_1{a}^{32}_1 + a^{24}_1{a}^{34}_1) = \\
   & e_{121}(e_1(e_2 + e_3) + (e_2 + e_4)(e_3 + e_4)) =  \\
   & e_1(e_3 + e_4(e_1 + e_2))(e_1(e_2 + e_3) + (e_2 + e_4)(e_3 + e_4))
    \subseteq e_{121}.
 \end{split}
\end{equation*}

  2) $\alpha = 321 = 341$,
  (Table \ref{table_adm_elem_D4}, Line $G31$, $t = 0, s = 0, r = 1$). We have
\begin{equation*}
 \begin{split}
   e_{321} = e_{341} = & e_3{a}^{21}_1{a}^{14}_1 = \\
    &  e_3(e_2 + e_1)(e_4 + e_1),  \\
  f_{3210} = f_{3410} = & e_{321}({a}^{14}_2 + e_3a^{24}_1) = \\
   & e_{321}(e_1 + e_4(e_3 + e_2) + e_3(e_2 + e_4)) =  \\
   & e_{321}(e_1 + (e_3 + e_2)(e_4 + e_2)(e_3 + e_4)) = \\
   & e_3(e_2 + e_1)(e_4 + e_1)(e_1 + (e_3 + e_2)(e_4 + e_2)(e_3 + e_4))
    \subseteq e_{321}.
 \end{split}
\end{equation*}

  \underline{For $n=4$}: \indent $\alpha = 2341 = 2321 = 2141$,
  (Table \ref{table_adm_elem_D4}, Line $F21$, $s = 0, t = 1, r = 1$). We have
\begin{equation*}
 \begin{split}
  e_{2141} = & e_2{a}^{41}_2{a}^{34}_1 = \\
     & e_2(e_4 + e_1(e_3 + e_2))(e_3 + e_4), \\
  f_{21410} = & e_{2141}(a^{34}_2 + e_1{a}^{24}_2) = \\
  & e_{2141}(e_2(e_4 + e_1(e_3 + e_2)) + e_1(e_2 + e_4(e_1 +  e_3)))
    \subseteq e_{21}.
 \end{split}
\end{equation*}

\subsection{Cumulative polynomials}
  \label{def_cumul}

According to Gelfand and Ponomarev we consider elements $e_t(n)$
and $f_0(n)$ used for construction of the perfect elements,
 \cite[p.7, p.53]{GP74}.

 The {\it cumulative  elements} $e_t(n)$, where $t = 1,2,3,4$ and
$f_0(n)$ are sums of all admissible elements of the same length
$n$, where $n$ is the length of the multi-index. So, for the case
$D^4$, the cumulative polynomials defined as follows:
\begin{equation}
 \label{cumul_polyn_D4}
 \begin{split}
   & e_t(n) = \sum{e}_{i_{n}\dots{i}_{2}t}, \qquad t = 1,2,3,4, \\
   & f_0(n) = \sum{f}_{i_{n}\dots{i}_{2}0},
 \end{split}
\end{equation}
where admissible elements $e_\alpha$ and $f_{\alpha0}$ are defined
in Table \ref{table_adm_elem_D4}, see also Section \ref{adm_D4},
Section \ref{examples_D4}.

{\it Examples of the cumulative polynomials in $D^4$}.

\underline{For $n=1$}:
\begin{equation*}
  \begin{split}
 & e_1(1) = e_1, \\
 & f_0(1) = I,
 \end{split}
\end{equation*}
  i.e., the {\it cumulative polynomials} for $n=1$ coincide with
the generators of $D^4$. \vspace{2mm}

\underline{For $n=2$}: \quad
  We use here the following property:
  \begin{equation*}
      e_\alpha({e_1{a}_t^{34} + {a}_s^{24}a_r^{32}}) =
      e_\alpha({a}_t^{43} + e_2{a}_{s-1}^{41}a_{r+1}^{31}),
  \end{equation*}
  see Lemma \ref{homom_polynom_P}, heading 2).
   For the type $G11$ in Table \ref{table_adm_elem_D4}, we have
  \begin{equation}
   \label{form_for_f_a0}
      e_\alpha({e_1{a}_{2t+1}^{34} + {a}_{2s+1}^{24}a_{2r-1}^{32}}) =
      e_\alpha({a}_{2t+1}^{43} + e_2{a}_{2s}^{41}a_{2r}^{31}),
  \end{equation}
  Taking in (\ref{form_for_f_a0}), $r = s = t = 0$ we have
  \begin{equation*}
      f_{10} =
      e_1({a}_{1}^{43} + e_2) = e_1(e_2 + e_3 + e_4).
  \end{equation*}
Thus,
\begin{equation}
 \label{cumul_ef2}
 \begin{split}
  e_1(2) = & e_{21} + e_{31} + e_{41} = \\
           & e_2(e_3 + e_4) + e_3(e_2 + e_4) + e_4(e_2 + e_3), \\
  f_0(2) = & f_{10} + f_{20} + f_{30} + f_{40} = \\
           & e_1(e_2 + e_3 + e_4) + e_2(e_1 + e_3 + e_4) +
             e_3(e_1 + e_2 + e_4) + e_4(e_1 + e_2 + e_3).
 \end{split}
\end{equation}

\underline{For $n=3$}:
\begin{equation*}
 \begin{split}
  e_1(3) = & e_{321} + e_{231} + e_{341} + e_{431} + e_{241} + e_{421}, \\
  & \\
  f_0(3) = & (f_{120} + f_{130} + f_{140}) +
             (f_{210} + f_{230} + f_{240}) + \\
           & (f_{310} + f_{320} + f_{340}) +
             (f_{310} + f_{320} + f_{340}),
 \end{split}
\end{equation*}
and so forth, see also Section \ref{examples_D4}.

\subsection{Perfect elements}

\subsubsection{Perfect cubes in the free modular lattice $D^r$}
  \label{cubes}
  \index{modular lattice!- $D^r$}
  \index{distributive lattice!- $B^+$}
  \index{distributive lattice!- $B^+(n)$}
  \index{Boolean cube(=perfect Boolean cube)}

In \cite{GP74}, Gelfand and Ponomarev constructed the sublattice
$B$ of perfect elements for the free modular lattice $D^r$ with
$r$ generators:
\begin{equation*}
   B = B^+ \bigcup B^+, \hspace{3mm} \text{ where }
   B^+ = \bigcup\limits_{n=1}^\infty B^+(n), \hspace{5mm}
   B^- = \bigcup\limits_{n=1}^\infty B^-(n).
\end{equation*}
They proved that every sublattice $B^+(n)$ (resp. $B^-(n)$) is
$2^r$-element Boolean algebra, so-called {\it Boolean cube} (which
can be also named {\it perfect Boolean cube}) and these {\it cubes
} are ordered in the following way. Every element of the cube
$B^+(n)$ is included in every element of the cube $B^+(n+1)$,
i.e.,
\begin{equation*}
 \left .
 \begin{array}{ll}
    v^+(n) \in B^+(n) \\
    v^+(n+1) \in B^+(n+1)
 \end{array}
 \right \}
 \Longrightarrow
 v^+(n+1) \subseteq v^+(n).
\end{equation*}
By analogy, the dual relation holds:
\begin{equation*}
 \left .
 \begin{array}{ll}
    v^-(n) \in B^-(n) \\
    v^-(n+1) \in B^-(n+1)
 \end{array}
 \right \}
 \Longrightarrow
 v^-(n) \subseteq v^-(n+1).
\end{equation*}

\subsubsection{Perfect elements in $D^4$}
  \label{sect_perfect_D4}
  Perfect elements
in $D^4$ (similarly, for $D^r$, see \cite[p.6]{GP74})
 are constructed by means of cumulative elements (\ref{cumul_polyn_D4})
 as follows:
\begin{equation}
  \label{perfect_D4}
   h_t(n) = \sum\limits_{j \neq t}{e_j(n)}.
\end{equation}
 Elements (\ref{perfect_D4}) generate perfect Boolean cube $B^+(n)$ from
 Section \ref{cubes}, see \cite{GP74}.

{\it Examples of perfect elements.}

 \underline{For $n=1$:}
\begin{equation}
  \label{perfect_D4_h1}
  \begin{split}
   & h_1(1) = e_2 + e_3 + e_4, \quad h_2(1) = e_1 + e_3 + e_4, \\
   & h_3(1) = e_1 + e_2 + e_4, \quad h_4(1) = e_1 + e_2 + e_3.
  \end{split}
\end{equation}

 \underline{For $n=2$}, by (\ref{cumul_ef2}) we have
\begin{equation}
  \label{perfect_D4_h1_2}
  \begin{split}
    h_1(2) = & e_2(2) + e_3(2) + e_4(2) = \\
    &(e_{12} + e_{32} + e_{42}) +
     (e_{13} + e_{23} + e_{43}) + (e_{14} + e_{24} + e_{34}) = \\
    & e_1(e_3 + e_4) + e_3(e_1 + e_4) + e_4(e_1 + e_3) + \\
    & e_1(e_2 + e_4) + e_2(e_1 + e_4) + e_4(e_1 + e_2) + \\
    & e_1(e_2 + e_3) + e_2(e_1 + e_3) + e_3(e_1 + e_2) = \\
    & (e_1 + e_3)(e_1 + e_4)(e_3 + e_4) +
      (e_1 + e_2)(e_1 + e_4)(e_2 + e_4) + \\
    &  (e_1 + e_2)(e_1 + e_3)(e_2 + e_3).
  \end{split}
\end{equation}

Let $h^{\max}(n)$ (resp. $h^{\min}(n)$) be the maximal (resp.
minimal) element in the cube $B^+(n)$.
\begin{equation}
  \label{perfect_D4_max_min}
  \begin{split}
    & h^{\max}(n) = \sum\limits_{i=1,2,3,4}h_i(n), \\
    & \\
    & h^{\min}(n) = \bigcap_{i=1,2,3,4}h_i(n). \\
  \end{split}
\end{equation}

\subsubsection{Perfect elements in $D^4$ by C.~Herrmann}
  C.~Herrmann introduced in \cite{H82}, \cite{H84}  polynomials $q_{ij}$ and associated
endomorphisms $\gamma_{ij}$ of $D^4$ playing the central role in
his construction of perfect polynomials.
 For $\{i,j,k,l\} = \{1, 2, 3, 4\}$, define
$$
   q_{ij} = q_{ji} = q_{kl} = q_{lk} = (e_i + e_j)(e_k + e_l),
$$
and
$$
   \gamma_{ij}f(e_1, e_2, e_3, e_4) = f(e_1q_{ij}, e_2q_{ij}, e_3q_{ij}, e_4q_{ij}),
$$
see Section \ref{sect_Herrmann_introduced}.

  C.~Herrmann's construction of perfect elements $s_n$, $t_n$ and
  $p_{i,n}$, where $i = 1,2,3,4$, is as follows:

  \begin{equation*}
   \label{def_stn_prelim}
   \begin{split}
      s_0 = I, & \quad s_1 = e_1 + e_2 + e_3 + e_4, \\
               & \quad s_{n+1} = \gamma_{12}(s_n) + \gamma_{13}(s_n) + \gamma_{14}(s_n), \\
      & \\
      t_0 = I, & \quad
         t_1 = (e_1 + e_2 + e_3)(e_1 + e_2 + e_4)(e_1 + e_3 + e_4)(e_2 + e_3 + e_4), \\
               & \quad t_{n+1} = \gamma_{12}(t_n) + \gamma_{13}(t_n) + \gamma_{14}(t_n), \\
      & \\
      p_{i,0} = I, & \quad
        p_{i,1} = e_i + t_1, \text{ where } i = 1,2,3,4, \\
      & \quad p_{i,n+1} = \gamma_{ij}(p_{j,n}) + \gamma_{ik}(p_{k,n}) + \gamma_{il}(p_{l,n}), \\
      & \quad
        \text{ where } i = 1,2,3,4, \text{ and } \{i,j,k,l\} = \{1,2,3,4\},
   \end{split}
  \end{equation*}
see Section \ref{sect_stp_poly}. Then,
    \begin{equation}
     \label{sn_as_sum}
      \begin{split}
       s_n =  &
        \sum\limits_{r + s + t = n-1}
        \gamma^t_{12}\gamma^s_{13}\gamma^r_{14}(\sum\limits_{i=1,2,3,4}e_i) =
         \sum\limits_{i=1,2,3,4}e_i(n) \simeq h^{\max}(n), \\
        & \\
       t_n = &
         \sum\limits_{r + s + t = n}
         \gamma^t_{12}\gamma^s_{13} \gamma^r_{14}
          (\sum\limits_{i=1,2,3,4}e_i(e_j + e_k + e_l)) =
            f_0(n+1) \simeq h^{\min}(n),  \\
        & \\
       p_{i,n} = & e_i(n) + f_0(n+1) \simeq h_j(n)h_k(n)h_l(n),
       \text{ for } \{i,j,k,l\} = \{1,2,3,4\},
      \end{split}
    \end{equation}
where $a \simeq b$ means modulo linear equivalence, see Theorem
\ref{prop_endom_R}, Theorem \ref{th_Herrmann_82}.
 Thus, the lattice of perfect elements generated by elements $s_n$, $t_n$,
$p_{i,n}$ introduced by Herrmann \cite{H82} coincide modulo linear
equivalence with the Gelfand-Ponomarev perfect elements, see Table
\ref{tab_perfect_in D4}.

 \begin{table} [h]
   \renewcommand{\arraystretch}{1.5}
  \begin{tabular} {||c|c|c|c||}
  \hline \hline
   N & Gelfand - Ponomarev     & Herrmann     & Sum of cumulative \\
     & definition              & definition   & polynomials \\
   \hline \hline
   1 & $\sum\limits_{i=1,2,3,4}h_i(n)$
     & $s_n$ & $\sum\limits_{i=1,2,3,4}e_i(n)$ \\
  \hline
   2 & $h_1(n)$ & $p_{2,n} + p_{3,n} + p_{4,n}$
     & $e_2(n) + e_3(n) + e_4(n)$ \\
  \hline
   3 & $h_2(n)$ & $p_{1,n} + p_{3,n} + p_{4,n}$
     & $e_1(n) + e_3(n) + e_4(n)$ \\
  \hline
   4 & $h_3(n)$ & $p_{1,n} + p_{2,n} + p_{4,n}$
     & $e_1(n) + e_2(n) + e_4(n)$  \\
  \hline
   5 & $h_4(n)$ & $p_{1,n} + p_{2,n} + p_{3,n}$
     & $e_1(n) + e_2(n) + e_3(n)$    \\
  \hline
   6 & $h_1(n)h_2(n)$
     & $p_{3,n} + p_{4,n}$  & $e_3(n) + e_4(n) + f_0(n+1)$   \\
  \hline
   7 & $h_1(n)h_3(n)$
     & $p_{2,n} + p_{4,n}$ &  $e_2(n) + e_4(n) + f_0(n+1)$    \\
  \hline
   8 & $h_1(n)h_4(n)$
     & $p_{2,n} + p_{3,n}$ &  $e_2(n) + e_3(n) + f_0(n+1)$    \\
  \hline
   9 & $h_2(n)h_3(n)$
     & $p_{1,n} + p_{4,n}$ &  $e_1(n) + e_4(n) + f_0(n+1)$    \\
  \hline
  10 & $h_2(n)h_4(n)$
     & $p_{1,n} + p_{3,n}$ &  $e_1(n) + e_3(n) + f_0(n+1)$    \\
  \hline
  11 & $h_3(n)h_4(n)$
  & $p_{1,n} + p_{2,n}$ &  $e_1(n) + e_2(n) + f_0(n+1)$    \\
  \hline
  12 & $h_1(n)h_2(n)h_3(n)$
     & $p_{4,n}$ &  $e_4(n) + f_0(n+1)$    \\
  \hline
  13 & $h_1(n)h_2(n)h_4(n)$
     & $p_{3,n}$ &  $e_3(n) + f_0(n+1)$    \\
  \hline
   14 & $h_1(n)h_3(n)h_4(n)$
      & $p_{2,n}$ &  $e_2(n) + f_0(n+1)$    \\
  \hline
   15 & $h_2(n)h_3(n)h_4(n)$
      & $p_{1,n}$ &  $e_1(n) + f_0(n+1)$    \\
  \hline
   16 & $\bigcap\limits_{i = 1,2,3,4}h_i(n)$
      & $t_n$ &  $f_0(n+1) = f^\vee_0(n+1)$    \\
  \hline \hline
  \end{tabular}
   \vspace{2mm}
   \caption{\hspace{3mm}Perfect elements in $B^+(n)$}
   \label{tab_perfect_in D4}
 \end{table}
  Table \ref{tab_perfect_in D4} gives two system of
  generators in the sublattice of the perfect elements $B^+(n)$ in $D^4$:
  coatoms $h_i(n)$ (by Gelfand-Ponomarev) vs. atoms $p_{i,n}$ (by
  C.~Herrmann).

%% file: 3admiss_D4.tex
\section{\sc\bf Atomic and admissible polynomials}
   \label{sect_adm_seq_D4}

\subsection{Admissible sequences}
  \label{subsect_adm_seq_D4}
 For definition of admissible sequences in the case of the modular
lattice $D^4$, see
 Section \ref{adm_D4}. Essentially, the fundamental property of this definition is
\begin{equation}
  \label{main_D4}
   ijk = ilk \indent \text{ for all } \{i,j,k,l\} = \{1,2,3,4\}.
\end{equation}
Relation (\ref{main_D4}) is our main tool in all further
calculations of admissible sequences of $D^4$.

 Without loss of generality only sequences starting at $1$ can
be considered. The following proposition will be used for the
classification of admissible sequences in $D^4$.

\begin{proposition}
 \label{relations_1_14}
  The following relations hold
 \item[1)]  $(31)^r(32)^s(31)^t = (32)^s(31)^{r+t},$
 \item[2)]  $(31)^r(21)^s(31)^t = (31)^{r+t}(21)^s,$
 \item[3)]  $(42)^r(41)^s = (41)^r(31)^s$, \quad $s \geq 1,$
 \item[4)]  $2(41)^r(31)^s = 2(31)^{s+1}(41)^{r-1},$
   \quad $r \geq 1,$
 \item[5)]  $(43)^r(42)^s(41)^t = (41)^r(21)^s(31)^t,$
 \item[6)]  $1(41)^r(21)^s = 1(21)^s(41)^r,$ \quad
            $1(i1)^r(j1)^s = 1(j1)^s(i1)^r,$ \quad
            $i,j \in \{2,3,4\}$,  $i \neq j,$
 \item[7)]  $(41)^r(21)^t(31)^s = (41)^r(31)^s(21)^t,$
 \item[8)]  $(13)^s(21)^r = (12)^r(31)^s,$
 \item[9)]  $12(41)^r(31)^s(21)^t =
   (14)^r(31)^{s+1}(21)^t = (14)^r(21)^{t+1}(31)^s,$
 \item[10)] $12(14)^r(31)^s(21)^t = (14)^r(31)^s(21)^{t+1},$
 \item[11)] $13(14)^r(31)^s(21)^t = (14)^r(31)^{s+2}(21)^{t-1},$
 \item[12)] $32(14)^r(31)^s(21)^t = (31)^s(21)^{t+1}(41)^r =
            34(14)^r(31)^s(21)^t,$
 \item[13)] $42(14)^r(31)^s(21)^t =
            (41)^s(21)^{t+1}(31)^r = 43(14)^r(31)^s(21)^t,$
 \item[14)] $23(14)^r(31)^s(21)^t = (21)^{t+1}(31)^s(41)^r =
            24(14)^r(31)^s(21)^t$,
 \item[15)] $2(41)^r(31)^s(21)^t = 2(14)^r(31)^{s+1}(21)^{t-1}$,
            \quad $t \geq 1, s > 1$.
\end{proposition}
For the proof, see \cite[Prop. 4.1.1]{St04}. \qedsymbol

\begin{table}[h]
 \renewcommand{\arraystretch}{1.9}
  \begin{tabular} {||c|c|c|c|c|c||}
  \hline \hline
     & Admissible Sequence
     & Action       & Action      & Action   & Action \cr
     &  & $\varphi_1$  & $\varphi_2$ & $\varphi_3$ & $\varphi_4$
     \\
  \hline         
    $F21$  & $(21)^t(41)^r(31)^s = (21)^t(31)^s(41)^r $
     & $G11$
     & -
     & $G31$
     & $G41$ \\
  \hline         
    $F31$  & $(31)^s(41)^r(21)^t = (31)^s(21)^t(41)^r $
     & $G11$
     & $G21$
     & -
     & $G41$ \\
  \hline      
    $F41$  & $(41)^r(31)^s(21)^t = (41)^r(21)^t(31)^s $
     & $G11$
     & $G21$
     & $G31$
     & - \\
  \hline \hline        
    $G11$  & $1(41)^r(31)^s(21)^t$ =
     & -
     & $F21$
     & $F31$
     & $F41$ \cr
     &    \hspace{3mm} $1(31)^s(41)^r(21)^t$ = & & & &  \cr
     &    \hspace{6mm} $1(21)^t(31)^s(41)^r$ & & & & \\
  \hline         
    $G21$  & $2(41)^r(31)^s(21)^t$ =
                $2(31)^{s+1}(41)^{r-1}(21)^t$ =
     & $H11$
     & -
     & $F31$
     & $F41$ \cr
     & \hspace{3mm} $2(14)^r(31)^{s+1}(21)^{t-1}$ & & & & \\
  \hline         
     $G31$  & $3(41)^r(21)^t(31)^s$ =
                $3(21)^{t+1}(41)^{r-1}(31)^s$ =
     & $H11$
     & $F21$
     & -
     & $F41$  \cr
     & \hspace{3mm} $3(14)^r(21)^{s+1}(31)^{t-1}$ & & & & \\
  \hline      
    $G41$  & $4(21)^t(31)^s(41)^r$ =
                $4(31)^{s+1}(21)^{t-1}(41)^r$ =
     & $H11$
     & $F21$
     & $F31$
     & -        \cr
     & \hspace{3mm} $4(12)^r(31)^{s+1}(41)^{t-1}$ & & & & \\
  \hline \hline     
    $H11$  & $(14)^r(31)^s(21)^t = (14)^r(21)^{t+1}(31)^{s-1}$ =
     & -
     & $G21$
     & $G31$
     & $G41$ \cr
     &   $(13)^s(41)^r(21)^t = (13)^s(21)^{t+1}(41)^{r-1}$ = & & & &  \cr
     &   $(12)^{t+1}(41)^r(31)^{s-1} = (12)^{t+1}(31)^{s-1}(41)^r$ & & & & \\
  \hline  \hline
  \end{tabular}
  \vspace{2mm}
  \caption{\hspace{3mm}Admissible sequences for the modular lattice $D^4$}
  \label{table_admissible_ExtD4}
\end{table}

 \begin{remark}[Note to Table \ref{table_admissible_ExtD4}]
  \label{def_Fij}
 {\rm 1)
Type $Fij$ (resp. $Gij$, $H11$) denotes the admissible sequence
starting at $j$ and ending at $i$. Sequences of type $Fij$ and
$H11$ contain an even number of symbols, sequences of type $Gij$
and $H11$ contain an odd number of symbols. For differences in
types $Fij, Gij, H11$, see the table.

 2) Thanks to heading 15) of Proposition \ref{relations_1_14}
 types $H21, H31, H41$ from \cite[p. 55, Table 4.1]{St04}
 are excluded, and the number of different cases of
 admissible sequences is equal to $8$ instead $11$ in \cite{St04}.
 \vspace{3mm}
 }
 \end{remark}

\begin{proposition}
 \label{full_adm_seq_D4}
   Full list of admissible sequences starting at $1$ is given
   by Table \ref{table_admissible_ExtD4}.
\end{proposition}

 \PerfProof It
suffices to prove that maps $\varphi_i$, where $i=1,2,3,4$, do not
lead out of Table
 \ref{table_admissible_ExtD4}. The exponents $r,s,t$ may be any
non-negative integer number. The proof is based on the relations
from Proposition \ref{relations_1_14}, see \cite[Prop.
4.1.2]{St04}.

\begin{figure}[h]
\centering
\includegraphics{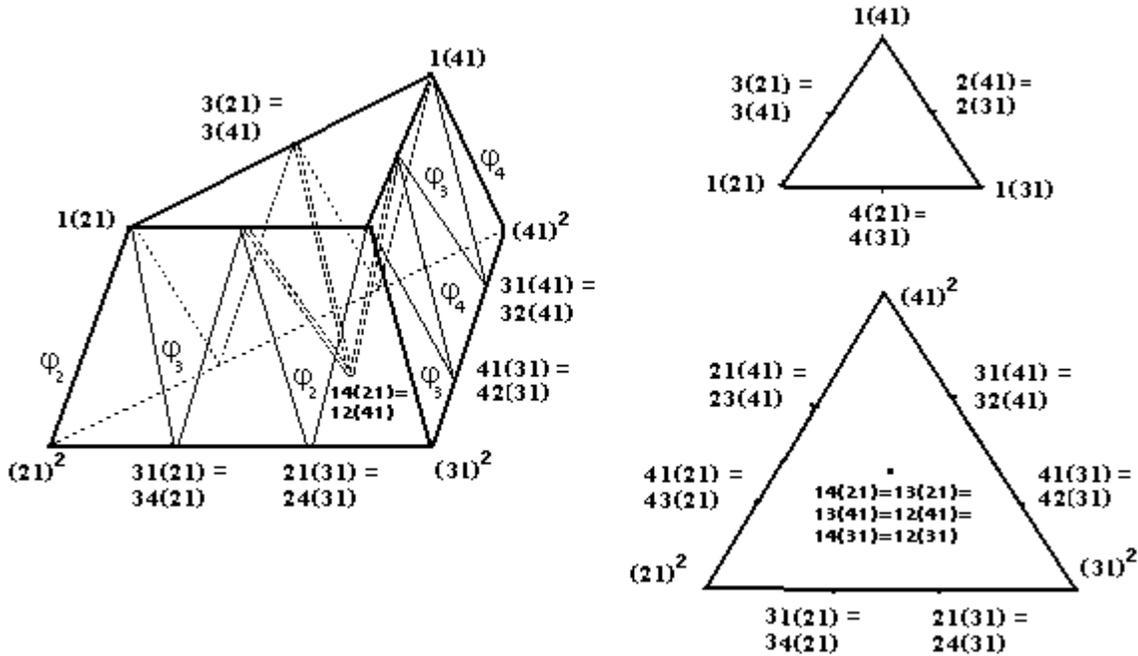}
\caption{\hspace{3mm}Slices of admissible sequences, $l=3$ and
$l=4$}
\label{section34}
\end{figure}

\begin{figure}[h]
\centering
\includegraphics{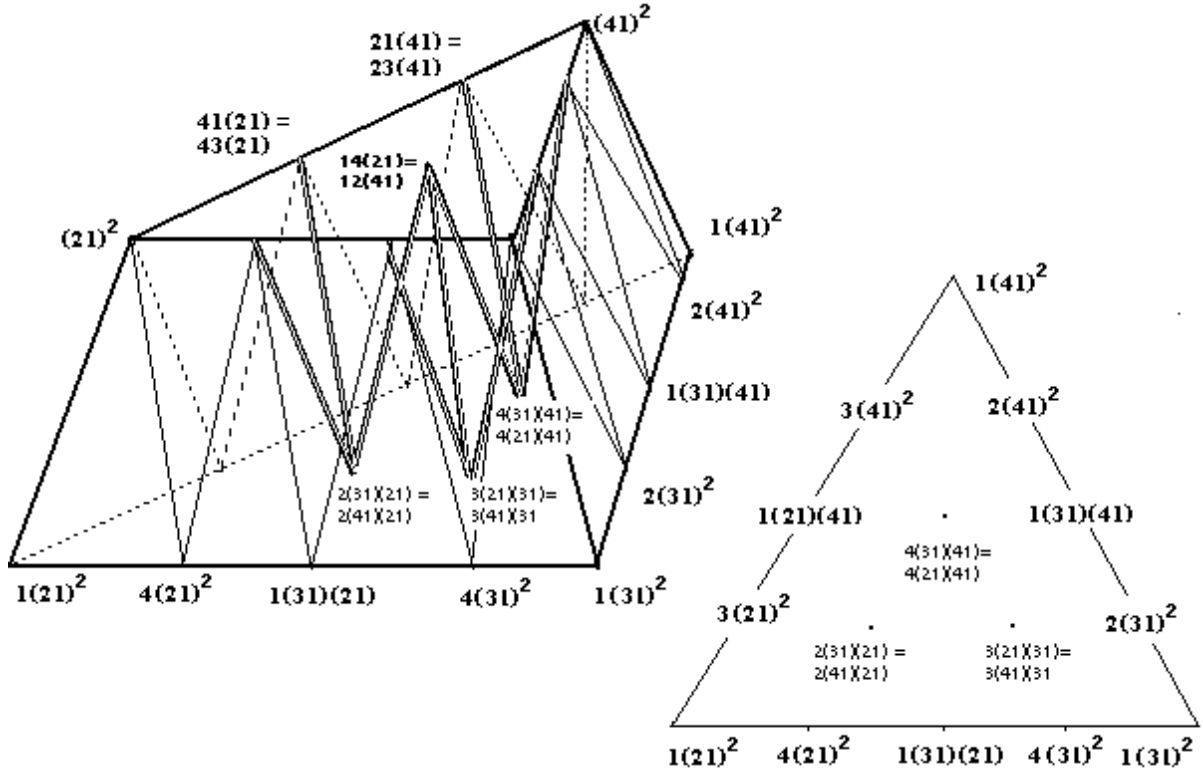}
\caption{\hspace{3mm}Slices of admissible sequences, $l=4$ and
$l=5$}
\label{section45}
\end{figure}

The pyramid on the Fig. \ref{pyramid_D4} has internal points. We
consider the slice $S(n)$ containing all sequences of the same
length $n$. The slices $S(3)$ and $S(4)$ are shown on the Fig.
\ref{section34}.
 The slice $S(4)$ contains only one internal point
\begin{equation}
  14(21) = 13(21) = 13(41) = 12(41) = 14(31) = 12(31).
\end{equation}
 The slices $S(4)$ and $S(5)$ are shown on the Fig. \ref{section45}.
 The slice $S(5)$ contains $3$ internal points
\begin{equation}
\begin{split}
  & 2(31)(21) = 2(41)(21), \\
  & 3(21)(31) = 3(41)(31), \\
  & 4(21)(41) = 4(31)(41).
\end{split}
\end{equation}

\begin{remark}
 \label{slices}
 {\rm The slice $S(n)$ contains $\displaystyle\frac{1}{2}{n(n+1)}$ different admissible
 sequences. Actions of $\varphi_i$, where $i = 1,2,3,4$ move
 every line in the triangle $S(n)$, which is parallel to
 some edge of the triangle, to the edge of $S(n+1)$. The edge
 containing $k$ points is moved to $k+1$ points in the $S(n+1)$. }
\end{remark}

\subsection{Atomic polynomials and elementary maps}
 \label{atomic_D4}

  \index{atomic elements $a^n_{ij}$ in $D^4$}

The free modular lattice $D^4$ is generated by $4$ generators:
$$
    D^4 = \{e_1,e_2,e_3,e_4\}.
$$

\begin{proposition}
  1) The following property of the atomic elements take place
   \begin{equation}
     \label{main_atomic_D4}
      e_j{a}_n^{kl} = e_j{a}_n^{lk} \text{ for } n \geq 1,
      \text{ and distinct indices } j,k,l.
   \end{equation}

  2) The definition of the atomic elements $a_n^{ij}$ in (\ref{def_atomic}) is
   well-defined. \vspace{2mm}

  3) We have
   \begin{equation}
     \label{incl_atomic_D4}
        a^{ij}_n \subseteq a^{ij}_{n-1} \subseteq \dots \subseteq
        a^{ij}_2 \subseteq a^{ij}_1 \subseteq a^{ij}_0 = I
        \text{ for all } i \neq j.
   \end{equation}

  4) To equalize the lower indices of the admissible polynomials
  $f_{{\alpha}0}$
  (see Table \ref{table_adm_elem_D4} and Theorem \ref{th_adm_classes_D4})
  we will use the following relation:
   \begin{equation}
     \label{equalize_D4}
        e_j + e_i{a}^{jk}_{t+1}{a}^{lk}_{s-1} =
        e_j + e_k{a}^{ij}_{t}{a}^{il}_{s}
        \quad \text{ for all } \{i,j,k,l\} = \{1,2,3,4\}.
   \end{equation}
\end{proposition}
  For the proof, see \cite[Prop. 4.2.1]{St04} \qedsymbol

Now we briefly recall definitions due to Gelfand and Ponomarev
\cite{GP74}, \cite{GP76}, \cite{GP77} of spaces $G_i$, $G_i{'}$,
representations $\nu^0$, $\nu^1$, joint maps $\psi_i$, and
elementary maps $\varphi_i$, where $i = 1,2,3,4$. To compare these
definitions with a case of the modular lattice $D^{2,2,2}$, see
 \cite[Sect. 2]{St04} \vspace{2mm}

We denote by
$$
  \{ Y_1, Y_2, Y_3, Y_4 \mid Y_i \subseteq X_0, i=1,2,3,4 \}
$$
the representation $\rho_X$ of $D^4$ in the finite dimensional
vector space
 $X = X_0$, and by
$$
  \{ Y^1_1, Y^1_2, Y^1_3, Y^1_4 \mid Y^1_i \subseteq X^1_0, i=1,2,3,4 \}
$$
the representation $\rho_{X^+}$ of $D^4$. Here $Y_i$ (resp.
$Y^1_i$) is the image of the generator $e_i$
 under the representation $\rho_X$ (resp. $\rho_{X^+}$).
\begin{equation}
\label{repr_rho_D4}
\begin{array}{ccccccccccccc}
   & & & Y_3 & & & &
   & & Y^1_3 & & \vspace{1mm}\\
   & & & \indownarrow & & & &
   & & \indownarrow & & \vspace{1mm} \\
   & Y_1 & \hspace{-1em}\hookrightarrow\hspace{-1em}
   & X_0 & \hspace{-1em}\hookleftarrow\hspace{-1em}
   & Y_2 & \hspace{11mm}
   & Y^1_1 & \hspace{-1em}\hookrightarrow\hspace{-1em}
   & X^1_0 & \hspace{-1em}\hookleftarrow\hspace{-1em}
   & Y^1_2 \vspace{1mm} \\
   & & & \inuparrow & & & &
   & & \inuparrow & & \vspace{1mm} \\
   & \rho_X      & & Y_4 & & & &
   \rho_{X^+}  & & Y^1_4 & & \\ \vspace{5mm}
\end{array}
\end{equation}

The space $X^+ = X^1_0$ --- the space of the representation
$\rho_{X^+}$
--- is
$$
   X^1_0 = \{ (\eta_1, \eta_2, \eta_3, \eta_4) \mid \eta_i \in Y_{i},
   \quad \sum\eta_i = 0 \},
$$
where $i \in \{1,2,3,4\}$. We set
$$
    R = \bigoplus\limits_\text{$i = 1,2,3,4$}Y_i,
$$
i.e.,
$$
  R = \{ (\eta_1, \eta_2, \eta_3, \eta_4) \mid \eta_i \in Y_i, \quad i = 1,2,3 \}.
$$
Then, $X^1_0 \subseteq R$.

 For the case of $D^4$, the spaces $G_i$ and $G{'}_i$
are introduced as follows:
\begin{equation}
  \begin{split}
       & G_1 = \{(\eta_1, 0, 0, 0) \mid \eta_1 \in Y_1\}, \quad
         G'_1 = \{(0, \eta_2, \eta_3, \eta_4) \mid \eta_i \in Y_i\}, \\
       & G_2 = \{(0, \eta_2, 0, 0) \mid \eta_2 \in Y_2\}, \quad
         G'_2 = \{(\eta_1, 0, \eta_3, \eta_4) \mid \eta_i \in Y_i\}, \\
       & G_3 = \{(0, 0, \eta_3, 0) \mid \eta_3 \in Y_3\},  \quad
         G'_3 = \{(\eta_1, \eta_2, 0, \eta_4) \mid \eta_i \in Y_i\}, \\
       & G_4 = \{(0, 0, 0, \eta_4) \mid \eta_4 \in Y_4\},  \quad
         G'_4 = \{(\eta_1, \eta_2, \eta_3, 0) \mid \eta_i \in Y_i\}. \\
  \end{split}
\end{equation}

For details, see \cite[p.43]{GP74}.

The associated representations $\nu_0$, $\nu_1$ in $R$ are defined
by Gelfand and Ponomarev \cite[eq.(7.2)]{GP74}:
\begin{equation}
 \begin{split}
   \label{assoc_repr_D4}
        & \nu^0(e_i) = X_0^{1} + G_i, \indent i = 1,2,3, \\
        & \nu^1(e_i) = X_0^{1}G'_i,  \indent i = 1,2,3.
  \end{split}
\end{equation}
Following \cite{GP74}, we introduce the elementary maps
$\varphi_i$:
 $$
 \varphi_i : X_0^{1}  \longrightarrow X_0,
  (\eta_1, \eta_2, \eta_3, \eta_4) \longmapsto \eta_i.
 $$
 From the definition we have
\begin{equation}  \label{sum_fi_D4}
        \varphi_1  + \varphi_2 + \varphi_3 + \varphi_4  = 0.
\end{equation}

\begin{table}[h]
  \renewcommand{\arraystretch}{1.7}
  \begin{tabular} {|| c | c | c ||}
   \hline \hline
       Notions & $D^{2,2,2}$ & $D^4$  \\
   \hline \hline 
       Generators
       & $\{x_1 \subseteq y_1,
            x_2 \subseteq y_2, x_3 \subseteq y_3\}$
       & $\{e_1, e_2, e_3, e_4\}$ \\
     \hline
       Atomic
       & $a_n^{ij} = x_i + y_j{a}_{n-1}^{jk}$,
       & $a_n^{ij} = e_i + e_j{a}^{kl}$,  \cr
       elements
       & $A_n^{ij} = y_i + x_j{A}_{n-1}^{ki},$
       & $\{i,j,k,l \} = \{1,2,3,4\}$ \cr
       & $\{i,j,k\} = \{1,2,3\}$ &  \\
     \hline
       Representation $\rho_X$
       & $\rho_X(x_i) = X_i, \rho_X(y_i) = Y_i$
       & $\rho_X(e_i) = Y_i$ \\
     \hline
       Space $X^1_0$
       & $X^1_0 =
        \{ (\eta_1, \eta_2, \eta_3) \mid \eta_i \in Y_{i} \},$
       & $X^1_0 =
          \{ (\eta_1, \eta_2, \eta_3, \eta_4) \mid \eta_i \in Y_{i} \},$ \cr
       & where $\eta_1 + \eta_2 + \eta_3 = 0$
       & where $\eta_1 + \eta_2 + \eta_3 + \eta_4 = 0$ \\
     \hline
       Representation $\rho_{X^+}$
       & $\rho_{X^+}(x_i) = X^1_i, \rho_{X^+}(y_i) = Y^1_i$
       & $\rho_{X^+}(y_i) = Y^1_i$ \\
     \hline
       Spaces $G_i$, $H_i$
       & $G_1 = \{(\eta_1, 0, 0) \mid \eta_1 \in Y_1\}$,
       & $G_1 = \{(\eta_1, 0, 0, 0) \mid \eta_1 \in Y_1\}$ \cr
       & $H_1 = \{(\xi_1, 0, 0) \mid \xi_1 \in X_1\}$,
       & $G_2 = \{(0, \eta_2, 0, 0) \mid \eta_2 \in Y_2\}$ \cr
       & $G_2 = \{(0, \eta_2, 0) \mid \eta_2 \in Y_2\}$,
       & $G_3 = \{(0, 0, \eta_3, 0) \mid \eta_3 \in Y_3\}$ \cr
       & $H_2 = \{(0, \xi_2, 0) \mid \xi_2 \in X_2\}$,
       & $G_4 = \{(0, 0, 0, \eta_4) \mid \eta_4 \in Y_4\}$ \cr
       & $G_3 = \{(0, 0, \eta_3) \mid \eta_3 \in Y_3\}$,
       & \cr
       & $H_3 = \{(0, 0, \xi_3) \mid \xi_3 \in X_3\}$,
       & \\
     \hline
       Joint maps $\psi_i$
       & $\psi_{i}(a) = X_0^{1} + G_{i}(H_{i}^{'} + \nu^{1}(a))$
       & $\psi_{i}(a) = X_0^{1} + G_{i}(G_{i}^{'} + \nu^{1}(a))$ \\
     \hline
       Quasi-
       & $\psi_i(a)\psi_i(b)$ =
       & $\psi_i(a)\psi_i(b)$ = \cr
       multiplicativity
       & $\psi_i((a + e_i)(b + x_j{x}_k))$
       & $\psi_i((a + e_i)(b + e_j{e}_k{e}_l)$ \\
     \hline
       Elementary maps $\varphi_i$
       & $\varphi_i : X_0^{1}  \longrightarrow X_0$,
       & $\varphi_i : X_0^{1}  \longrightarrow X_0$, \cr
       & $(\eta_1, \eta_2, \eta_3) \longmapsto \eta_i$
       & $(\eta_1, \eta_2, \eta_3, \eta_4) \longmapsto \eta_i$ \\
     \hline
       Fundamental properties
       & $\varphi_i\varphi_j\varphi_i + \varphi_i\varphi_k\varphi_i = 0$,
       & $\varphi_i\varphi_k\varphi_j + \varphi_i\varphi_l\varphi_j = 0$ \cr
       of elementary maps $\varphi_i$
       & $\varphi_i^3 = 0$
       & $\varphi_i^2 = 0$ \\
     \hline
       Fundamental properties & &  \cr
       of indices
       & \fbox{$iji = iki$}
       & \fbox{$ikj = ilj$} \cr
       (admissible sequences) && \\
     \hline \hline
  \end{tabular}
  \vspace{2mm}
\caption{\hspace{3mm}Comparison of notions in $D^{2,2,2}$ and
$D^4$}
  \label{compare_notions}
  \vspace{3.7mm} 
\end{table}

\subsection{Basic relations for elementary and joint maps}
 \label{sect_basic_rel_D4}

 We define {\it joint} maps $\psi_i\colon{D}^4 \longrightarrow
\mathcal{L}(R)$ as in the case of $D^{2,2,2}$:
\begin{equation}
  \label{psi_D4}
    \psi_{i}(a) = X_0^{1} + G_i(G'_i + \nu^1(a)).
\end{equation}

\begin{proposition}
 \label{basic_eq_D4}
In the case $D^4$ the joint maps $\psi_i$ satisfy the following
basic relations:
\begin{enumerate}
 \item $\psi_i(e_i) = X_0^1$, \vspace{2mm}
 \item $\psi_i(e_j) = \nu^0(e_i(e_k + e_l))$, \vspace{2mm}
 \item $\psi_i(I) = \nu^0(e_i(e_j + e_k + e_l))$, \vspace{2mm}
 \item $\psi_i(e_k{e}_l) = \psi_j(e_k{e}_l) =
     \nu^0(e_i{e}_j)$. \vspace{2mm}
\end{enumerate}
\end{proposition}

\PerfProof 1) From (\ref{psi_D4}) and (\ref{assoc_repr_D4}) we
have $\psi_i(e_i) = X_0^1 + G_i{G'_i} = X_0^1$. \vspace{2mm}

2) We have
   $\psi_i(y_j) = X_0^1 + G_i(G'_i + X_0^1{G'_j})$.
       From $G_i \subseteq G'_j$ for $i \neq j$ and by
       the permutation property \cite[Sect. 2]{St04}, we get
$$
       \psi_i(y_j) = X_0^1 + G_i(G'_i{G'_j} + X_0^1).
$$
       Since $G'_iG'_j =  G_k + G_l$, where $i,j,k,l$ are distinct indices,
       we have
\begin{equation*}
 \begin{split}
       \psi_i(y_j) = &  X_0^1 + G_i(X_0^1 + G_k + G_l) = \\
       & (X_0^1 + G_i)(X_0^1 + G_k + G_l) = \nu^0(e_i(e_k + e_l)).
        \vspace{2mm}
 \end{split}
\end{equation*}

3) Again,
\begin{equation}
 \begin{split}
     \psi_i(I) = & X_0^1 + G_i(G'_i + X_0^1) = \\
     & X_0^1 + G_i(G_j + G_k + G_l + X_0^1) =
      \nu^0(e_i(e_j + e_k + e_l)).
     \vspace{2mm}
 \end{split}
\end{equation}

4) Since $G_i \subseteq G'_k$ for all $i \neq k$, we have
 \begin{equation}
   \label{e_k__e_l}
 \begin{split}
     \psi_i(e_k{e}_l) = & X_0^1 + G_i(G'_i + X_0^1G'_kG'_l) = \\
     & X_0^1 + G_i(G'_i{G'_k}G'_l + X_0^1).
 \end{split}
\end{equation}
Since $G_j(G_i + G_k + G_l) = 0$, then
\begin{equation}
   \label{intersect_G_2}
    G'_i{G'_j} = G_k + G_l.
\end{equation}
Indeed,
 \begin{equation}
   \label{intersect_G_3}
 \begin{split}
     G'_i{G'_j} = & (G_j + G_k + G_l)(G_i + G_k + G_l) = \\
     & G_k + G_l +  G_j(G_i + G_k + G_l) = G_k + G_l.
 \end{split}
\end{equation}
From (\ref{e_k__e_l}) and (\ref{intersect_G_2}), we see that
 \begin{equation*}
 \begin{split}
     \psi_i(e_k{e}_l) = & X_0^1 + G_i(G'_i{G'_k}G'_l + X_0^1) = \\
     & X_0^1 + G_i(X_0^1 + G_j) = \nu_0(e_i{e}_j).
     \qed \vspace{2mm}
 \end{split}
\end{equation*}

The main relation between the elementary map $\varphi_i$ and the
joint map $\psi_i$ (see \cite[Prop. 2.4.3]{St04}) holds also for
the case of $D^4$. Namely, let $a,b,c \subseteq D^4$, then

\begin{equation}
  \label{phi_and_psi_D4}
\begin{split}
    & \text{ (i) \qquad If } \psi_i(a) = \nu^0(b), \text{ then }
      \varphi_i\rho_{X^+}(a) = \rho_X(b). \\
    & \text{(ii) \qquad If } \psi_i(a) = \nu^0(b) \text{ and }
      \psi_i(ac) = \psi_i(a)\psi_i(c), \text{ then }\\
    & \qquad \qquad
      \varphi_i\rho_{X^+}(ac) = \varphi_i\rho_{X^+}(a)\varphi_i\rho_{X^+}(c).
 \end{split}
\end{equation}

From Proposition \ref{basic_eq_D4} and eq.(\ref{phi_and_psi_D4})
we have

\begin{corollary}
 \label{cor_psi_D4}
 For the elementary map $\varphi_i$ the following basic
 relations hold:
\begin{enumerate}
 \item $\varphi_i\rho_{X^+}(e_i) = 0$, \vspace{2mm}
 \item $\varphi_i\rho_{X^+}(e_j) = \rho_X(e_i(e_k + e_l))$, \vspace{2mm}
 \item $\varphi_i\rho_{X^+}(I) = \rho_X(e_i(e_j + e_k + e_l))$, \vspace{2mm}
 \item $\varphi_i\rho_{X^+}(e_k{y}_l) = \rho_X(e_i{e}_j)$. \vspace{2mm}
 \end{enumerate}
  \qedsymbol
\end{corollary}

\subsection{Additivity and multiplicativity of the joint maps}
 \label{add_multi_D4}

\begin{proposition}
  \label{additivity_D4}
   The map $\psi_i$ is additive and quasimultiplicative
   with respect to the
   lattice operations + and $\cap$,
   namely:
\begin{enumerate}
   \item $\psi_i(a) + \psi_i(b) = \psi_i(a+b)$, \vspace{2mm}
   \item $\psi_i(a)\psi_i(b)$ =
      $\psi_i((a + e_i)(b + x_j{x}_k{x}_l))$, \vspace{2mm}
   \item $\psi_i(a)\psi_i(b)$ =
      $\psi_i(a(b + e_i + x_j{x}_k{x}_l))$. \vspace{2mm}
\end{enumerate}
\end{proposition}
  For the proof, see \cite[Prop. 4.4.1]{St04} \qedsymbol

We need the following corollary (atomic multiplicativity) from
Proposition \ref{additivity_D4}
\begin{corollary}
  \label{cor_mul_d4}
  {\rm(a)} Suppose one of the following
  inclusions holds:
\begin{equation*}
\begin{array}{lll}
   & {\rm(i)} & \quad e_i + e_j{e}_k{e}_l \subseteq a, \\
   & {\rm (ii)} & \quad e_i + e_j{e}_k{e}_l \subseteq b,  \\
   & {\rm (iii)} & \quad e_i \subseteq a, \quad
              e_j{e}_k{e}_l \subseteq b, \\
   & {\rm (iv)} & \quad e_i \subseteq b, \quad
              e_j{e}_k{e}_l \subseteq a. \\
\end{array}
\end{equation*}
 Then the joint map $\psi_i$ operates as a homomorphism on the elements
 $a$ and $b$ with respect to the
 lattice operations $+$ and $\cap$, i.e.,
$$
   \psi_i(a) + \psi_i(b) = \psi_i(a+b), \hspace{3mm}
    \psi_i(a)\psi_i(b) = \psi_i(a)\psi_i(b).
$$

{\rm(b)} The joint map ${\psi_i}$ applied to the following atomic
elements
    is the intersection preserving map, i.e., multiplicative with respect
    to the operation $\cap$:
\begin{equation}
     \label{homo1_D4}
     \psi_i(ba_n^{ij}) = \psi_i(b)\psi_i(a_n^{ij})
             \text{ for every } b \subseteq D^4.
\end{equation}
\end{corollary}

\subsection{The action of maps $\psi_i$ and $\varphi_i$
            on the atomic elements}
\begin{proposition}
   \label{action_psi_D4}
  The joint maps $\psi_i$ applied to the atomic elements $a^{ij}_n$
 satisfy the
 following relations
 \begin{enumerate}
   \item $\psi_i(a^{ij}_n) = \nu^0(e_i{a}^{kl}_n)$, \vspace{2mm}
   \item $\psi_j(a^{ij}_n) = \nu^0(e_j(e_k + e_l))$, \vspace{2mm}
   \item $\psi_j(e_i{a}^{kl}_n) = \nu^0(e_j{a}^{kl}_{n+1})$. \vspace{2mm}
  \end{enumerate}
\end{proposition}

  For the proof, see \cite[Prop. 4.5.1]{St04}

\subsection{The fundamental property of the elementary maps}
 \label{sect_fundam_prop_D4}

\begin{proposition}
 \label{motiv_admis_D4}
  For $\{i,j,k,l\} = \{1,2,3,4\}$ the following
  relations hold
\begin{equation}
  \label{fund_varphi_rel_D4}
     \varphi_i\varphi_k\varphi_j +  \varphi_i\varphi_l\varphi_j = 0, \vspace{2mm}
\end{equation}
\begin{equation}
     \varphi_i^2 = 0.
\end{equation}
\end{proposition}
\PerfProof
  For every vector $v \in X_0^1$, by definition of $\varphi_i$
  we have $(\varphi_i + \varphi_j + \varphi_k + \varphi_l)(v) = 0$,
  see eq. (\ref{sum_fi_D4}).
  In other words, $\varphi_i + \varphi_j + \varphi_k + \varphi_l$ = 0. Therefore,
$$
    \varphi_i\varphi_k\varphi_j +  \varphi_i\varphi_l\varphi_j =
    \varphi_i(\varphi_i + \varphi_j)\varphi_j =
    \varphi_i^2\varphi_j + \varphi_i\varphi_j^2.
$$
  So, it suffices to prove that $\varphi_i^2$ = 0.
  For every
   $z \subseteq D^4$, by Corollary \ref{cor_psi_D4}, headings (3) and (1), we have
\begin{equation*}
\begin{split}
  & \varphi_i^2\rho_{X^2}(z) \subseteq
      \varphi_i(\varphi_i\rho_{X^2}(I)) =  \\
  & \varphi_i(\rho_{X^+}(e_i(e_j + e_k + e_l)))
      \subseteq \varphi_i\rho_{X^+}(e_i) = 0.
  \qed \vspace{2mm}
\end{split}
\end{equation*}

 \begin{corollary}
   \label{basic_rel_D4}
     The relation
 \begin{equation}
   \label{fund_D4}
      \varphi_i\varphi_k\varphi_j(B)
       = \varphi_i\varphi_l\varphi_j(B)
 \end{equation}
 takes place for every subspace $B \subseteq X_0^2$,
 where $X^2 = X_0^2$ is the representation space
 of $\rho_{X^2}$.
 \end{corollary}

Essentially, relations (\ref{fund_varphi_rel_D4}) and
(\ref{fund_D4}) are fundamental and motivate the construction of
the admissible sequences satisfying the following relation:
\begin{equation}
  \label{fund_indices}
   ikj = ilj,
\end{equation}
where indices $i,j,k,l$ are all distinct, see Table
 \ref{compare_notions} and Section \ref{subsect_adm_seq_D4}.

\subsection{The $\varphi_i-$homomorphic elements}
  \label{subs_homom_D4}
By analogy with the modular lattice $D^{2,2,2}$
 (see \cite[prop. Sect. 4.7]{St04}),
 we introduce now $\varphi_i-$homomorphic polynomials in $D^4$.

\index{$\varphi_i-$homomorphic polynomial}

   An element $a \subseteq D^4$ is said to be
    {\it $\varphi_i-$homomorphic}, if
    \begin{equation}
       \varphi_i\rho_{X^+}(ap) =
          \varphi_i\rho_{X^+}(a)\varphi_i\rho_{X^+}(p)
            \text{ for all } p \subseteq D^4.
    \end{equation}

    An element $a \subseteq D^4$ is said to be
    {\it $(\varphi_i, e_k)-$homomorphic}, if
    \begin{equation}
        \varphi_i\rho_{X^+}(ap) =
           \varphi_i\rho_{X^+}(e_k{a})\varphi_i\rho_{X^+}(p)
            \text{ for all } p \subseteq e_k.
    \end{equation}

\begin{theorem}
  \label{th_homomorhism_D4}
  1) The polynomials $a^{ij}_n$
  are $\varphi_i-$homomorphic.
    \vspace{2mm}

  2) The polynomials $a^{ij}_n$ are ($\varphi_j,e_k)-$homomorphic
     for distinct indices $\{i,j,k\}$.
            \vspace{2mm}
\end{theorem}
  For the proof, see \cite[Th. 1.7.1]{St04}.  \qedsymbol

\subsection{The theorem on the classes of admissible elements}
 \label{sect_adm_classes_D4}
\begin{theorem}
   \label{th_adm_classes_D4}
   Let $\alpha = i_n{i}_{n-1}\dots{1}$
   be an admissible sequence for $D^4$ and $i \neq i_n$.
   Then ${i}\alpha$ is admissible and, for
    $z_\alpha = e_\alpha$ or $f_{\alpha0}$ from
    Table \ref{table_adm_elem_D4},
    the following relation holds:
 \begin{equation}
   \label{adm_classes_D4}
    \varphi_i\rho_{X^+}(z_\alpha) = \rho_X(z_{i\alpha}).
 \end{equation}
\end{theorem}

For the proof of the theorem on admissible elements in $D^4$, see
 Section B.2 in \cite{St04}.

 The proof repeatedly uses the basic properties of the
admissible sequences in $D^4$ considered in Section
\ref{basic_adm_D4}, Lemma \ref{homom_polynom_P}.

\begin{table}[h]
 \renewcommand{\arraystretch}{1.3}
  \begin{tabular} {||c|c|c|c||}
  \hline \hline
     & Admissible
     & Admissible & Admissible  \cr
     & sequence $\alpha$ & polynomial $e_\alpha$
     & polynomial $f_{{\alpha}0}$ \\
  \hline         
    $F21$  & $(21)^t(41)^r(31)^s = $
     & $e_2{a}^{31}_{2s}a^{41}_{2r}a^{34}_{2t-1}$
     & $e_\alpha(e_2{a}^{34}_{2t} + a^{41}_{2r+1}a^{31}_{2s-1}) = $ \cr
     & $(21)^t(31)^s(41)^r $ &
     & $e_\alpha({a}^{43}_{2t} + e_1a^{24}_{2r}a^{23}_{2s})$ \\
  \hline         
    $F31$  & $(31)^s(41)^r(21)^t = $
     & $e_3{a}^{21}_{2t}a^{41}_{2r}a^{24}_{2s-1}$
     & $e_\alpha(e_3{a}^{42}_{2s} + a^{41}_{2r+1}a^{21}_{2t-1}) = $ \cr
     & $(31)^s(21)^t(41)^r  $ &
     & $e_\alpha({a}^{42}_{2s} + e_1a^{34}_{2r}a^{32}_{2t})$ \\
  \hline      
    $F41$  & $(41)^r(31)^s(21)^t = $
     & $e_4{a}^{21}_{2t}a^{31}_{2s}a^{32}_{2r-1}$
     & $e_\alpha(e_4{a}^{32}_{2r} + a^{31}_{2s+1}a^{21}_{2t-1}) = $ \cr
     & $(41)^r(21)^t(31)^s$ &
     & $e_\alpha({a}^{32}_{2r} + e_1a^{43}_{2s}a^{42}_{2t})$ \\
  \hline \hline        
    $G11$  & $1(41)^r(31)^s(21)^t$ = & &  \cr
     &   \hspace{3mm} $1(31)^s(41)^r(21)^t$ =
     & $e_1{a}^{24}_{2s}a^{34}_{2t}a^{32}_{2r}$
     & $e_\alpha(e_1{a}^{32}_{2r+1} + a^{24}_{2s+1}a^{34}_{2t-1}) = $ \cr
     &   \hspace{6mm} $1(21)^t(31)^s(41)^r$ &
     & $e_\alpha({a}^{32}_{2r+1} + e_4a^{21}_{2s}a^{31}_{2t})$   \\
  \hline         
    $G21$  & $2(41)^r(31)^s(21)^t =$
     & $e_2{a}^{34}_{2t}a^{31}_{2s+1}a^{14}_{2r-1} =$
     & $e_\alpha(e_2{a}^{14}_{2r} + a^{31}_{2s+2}a^{34}_{2t-1}) =$ \cr
     & $2(31)^{s+1}(41)^{r-1}(21)^t$
     & $e_2{a}^{34}_{2t}a^{31}_{2s-1}a^{14}_{2r+1}$
     & $e_\alpha({a}^{14}_{2r} + e_3a^{21}_{2s+1}a^{24}_{2t})$  \\
  \hline         
     $G31$  & $3(41)^r(21)^t(31)^s =$
     & $e_3{a}^{24}_{2s}a^{21}_{2t+1}a^{14}_{2r-1} =$
     & $e_\alpha(e_3{a}^{14}_{2r} + a^{21}_{2t+2}a^{24}_{2s-1}) =$ \cr
     & $3(21)^{t+1}(41)^{r-1}(31)^s$
     & $e_3{a}^{24}_{2s}a^{21}_{2t-1}a^{14}_{2r+1}$
     & $e_\alpha({a}^{14}_{2r} + e_2a^{31}_{2t+1}a^{34}_{2s})$ \\
  \hline      
    $G41$  & $4(21)^t(31)^s(41)^r =$
     & $e_4{a}^{32}_{2r}a^{31}_{2s+1}a^{12}_{2t-1} =$
     & $e_\alpha(e_4{a}^{12}_{2t} + a^{31}_{2s}a^{32}_{2r+1}) = $ \cr
     & $4(31)^{s+1}(21)^{t-1}(41)^r$
     & $e_4{a}^{32}_{2r}a^{31}_{2s-1}a^{12}_{2t+1}$
     & $e_\alpha({a}^{12}_{2t} + e_3a^{41}_{2s+1}a^{42}_{2r})$ \\
  \hline \hline     
    $H11$  & $(14)^r(31)^s(21)^t$ =
     & & \cr
     & $(14)^r(21)^{t+1}(31)^{s-1}$ =
     & $e_1a^{23}_{2r+1}{a}^{24}_{2s-1}a^{34}_{2t-1} = $
     & $e_\alpha(e_1a^{23}_{2r} + a^{24}_{2s}a^{34}_{2t}) = $ \cr
     & $(13)^s(41)^r(21)^t =$
     & $e_1a^{23}_{2r-1}{a}^{24}_{2s-1}a^{34}_{2t+1} = $
     & $e_\alpha(e_1a^{24}_{2s} + a^{34}_{2t}a^{23}_{2r}) = $ \cr
     & $(13)^s(21)^{t+1}(41)^{r-1}$
     & $e_1a^{23}_{2r-1}{a}^{24}_{2s+1}a^{34}_{2t-1}$
     & $e_\alpha(e_1a^{34}_{2t} + a^{23}_{2r}a^{24}_{2s})$ \cr
     & $(12)^{t+1}(41)^r(31)^{s-1}$ = & & \cr
     & $(12)^{t+1}(31)^{s-1}(41)^r$ & & \\
  \hline  \hline
  \end{tabular}
  \vspace{2mm}
  \caption{\hspace{3mm}Admissible polynomials in the modular lattice $D^4$}
  \footnotesize
  \begin{tabular}{l}
    Notes to Table: \cr
  1) For more details about admissible sequences,
     see Proposition \ref{full_adm_seq_D4} and
     Table \ref{table_admissible_ExtD4}. \cr
  2) For relations given in two last columns (definitions of admissible
     polynomials $e_\alpha$ and $f_{{\alpha}0}$), \cr
     see Lemma \ref{homom_polynom_P}. \cr
  3) In each line, each low index should be non-negative. For
     example, for Line $F21$, \cr
     we have:
     $s \geq 0, r \geq 0, t \geq 1$;
     for Line $G21$, we have: $s \geq 0, r \geq 0, t \geq 0$.
  \label{table_adm_elem_D4}
  \end{tabular}
 \normalsize
\end{table}

\begin{table}[h]
 \renewcommand{\arraystretch}{1.3}
  \begin{tabular} {|c|c|c|}
  \hline \hline
     N
     & \quad Equivalent form of $f_{{\alpha}0}$ \quad
     & Forms obtained by $e_j{a}^{kl}_n = e_j{a}^{lk}_n$ \\
  \hline         
      1
      & $e_\alpha(e_2{a}^{34}_{2t} + a^{41}_{2r+1}a^{31}_{2s-1})$
      & $e_\alpha(e_2{a}^{43}_{2t} + a^{41}_{2r+1}a^{31}_{2s-1})$ \\
  \hline
      2
      & $e_\alpha({a}^{34}_{2t} + e_2a^{41}_{2r+1}a^{31}_{2s-1})$
      & $e_\alpha({a}^{43}_{2t} + e_2a^{41}_{2r+1}a^{31}_{2s-1})$ \\
  \hline
      3
      & $e_\alpha({a}^{34}_{2t} + e_2a^{14}_{2r+1}a^{13}_{2s-1})$
      & $e_\alpha({a}^{43}_{2t} + e_2a^{14}_{2r+1}a^{13}_{2s-1})$ \\
  \hline
      4
      & $e_\alpha(e_2{a}^{34}_{2t} + a^{14}_{2r+1}a^{13}_{2s-1})$
      & $e_\alpha(e_2{a}^{43}_{2t} + a^{14}_{2r+1}a^{13}_{2s-1})$ \\
  \hline
      5
      & $e_\alpha(a^{14}_{2r+1} + e_2{a}^{34}_{2t}a^{13}_{2s-1})$
      & $e_\alpha(a^{14}_{2r+1} + e_2{a}^{43}_{2t}a^{13}_{2s-1})$ \\
  \hline
      6
      & $e_\alpha(e_2a^{14}_{2r+1} + {a}^{34}_{2t}a^{13}_{2s-1})$
      & $e_\alpha(e_2a^{14}_{2r+1} + {a}^{43}_{2t}a^{13}_{2s-1})$ \\
  \hline
      7
      & $e_\alpha({a}^{34}_{2t} + e_2a^{14}_{2r-1}a^{13}_{2s+1})$
      & $e_\alpha({a}^{43}_{2t} + e_2a^{14}_{2r-1}a^{13}_{2s+1})$ \\
  \hline
      8
      & $e_\alpha({a}^{34}_{2t}a^{14}_{2r-1} + e_2a^{13}_{2s+1})$
      & $e_\alpha({a}^{43}_{2t}a^{14}_{2r-1} + e_2a^{13}_{2s+1})$ \\
  \hline
      9
      & $e_\alpha(a^{13}_{2s+1} + e_2{a}^{34}_{2t}a^{14}_{2r-1})$
      & $e_\alpha(a^{13}_{2s+1} + e_2{a}^{43}_{2t}a^{14}_{2r-1})$ \\
  \hline  \hline
  \end{tabular}
  \vspace{2mm}
  \caption{\hspace{3mm}Different equivalent forms of the element
    $f_{{\alpha}0} = f_{(21)^t(41)^r(31)^s0}$}
  \label{equival_forms_D4}
\end{table}

\subsubsection{Basic properties of admissible elements in $D^4$}
  \label{basic_adm_D4}
We prove here a number of basic properties. of the atomic elements
in $D^4$ used in the proof of the theorem on admissible elements
(Theorem \ref{th_adm_classes_D4}). In particular, in some cases
the lower indices of polynomials $a^{ij}_s$ entering in the
admissible elements $f_{{\alpha}0}$ can be transformed as in the
following
\begin{lemma}
 \label{homom_polynom_P}
 1) Every polynomial $f_{{\alpha}0}$ from Table
 \ref{table_adm_elem_D4} can be represented as an intersection of
$e_\alpha$ and $P$. For every $i \neq i_n$
 (see Section \ref{adm_D4}),
 we select $P$ to be some $\varphi_i-$homomorphic
polynomial.

2) The lower indices of polynomials $a^{ij}_s$ entering in the
admissible elements $f_{{\alpha}0}$ can be equalized as follows:
\begin{equation}
 \label{sym_forms_P}
  \begin{split}
      f_{{\alpha}0} =
    & f_{(21)^t(41)^r(31)^s0} = \\
    & e_\alpha(e_2{a}^{34}_{2t} + a^{41}_{2r+1}a^{31}_{2s-1}) =
      e_\alpha({a}^{43}_{2t} + e_1a^{24}_{2r}a^{23}_{2s}).
  \end{split}
\end{equation}
The generic relation\footnote{Throughout this lemma
    we suppose that $\{i,j,k,l\} = \{1,2,3,4\}$.}
 is the following:
\begin{equation}
 \label{sym_forms_P_2}
  \begin{split}
    e_i(e_i{a}^{jl}_{t} + & a^{kj}_{r+1}a^{kl}_{s-1}) = \\
    & e_i({a}^{jl}_{t} + e_i{a}^{kj}_{r+1}a^{kl}_{s-1}) =
      e_i({a}^{jl}_{t} + e_k{a}^{ij}_{r}a^{il}_{s}).
  \end{split}
\end{equation}

3) The substitution
\begin{equation}
 \label{subst_r_s}
  \begin{cases}
     r \mapsto r-2, \\
     s \mapsto s+2
  \end{cases}
\end{equation}
does not change the polynomial $e_i(e_i{a}^{jl}_{t} +
a^{kj}_{r+1}a^{kl}_{s-1})$, namely:
\begin{equation}
 \label{sym_forms_P_3}
    e_i(e_i{a}^{jl}_{t} + a^{kj}_{r+1}a^{kl}_{s-1}) =
    e_i(e_i{a}^{jl}_{t} + a^{kj}_{r-1}a^{kl}_{s+1}).
\end{equation}

4) The substitution (\ref{subst_r_s}) does not change the
polynomial $e_ia^{kj}_{s}a^{kl}_{r}$:
\begin{equation}
 \label{sym_forms_P_4}
    e_ia^{kj}_{r+1}a^{kl}_{s-1} =
    e_ia^{kj}_{r-1}a^{kl}_{s+1}.
\end{equation}

\end{lemma}
  For the proof
  of this lemma see \cite[Lemma 4.8.2]{St04}.
 \qedsymbol \vspace{2mm}

\subsubsection{Coincidence with the Gelfand-Ponomarev polynomials in $D^4/\theta$}
  \label{coinc_GP_D4}
Let $e_\alpha, f_{{\alpha}0}$ be admissible elements constructed
in this work and  $\tilde{e}_\alpha, \tilde{f}_{{\alpha}0}$ be the
admissible elements constructed by Gelfand and Ponomarev
\cite{GP74}. We will prove that for the admissible sequences of
the small length, the coincidence of $e_\alpha$ with
$\tilde{e}_\alpha$ (resp. $f_{{\alpha}0}$ with
$\tilde{f}_{{\alpha}0}$) takes place in $D^4$ (not only in
$D^4/\theta$). Since Theorem \ref{th_adm_classes_D4} takes place
for both $e_\alpha, f_{{\alpha}0}$ and for $\tilde{e}_\alpha,
\tilde{f}_{{\alpha}0}$ \cite[Th.7.2, Th.7.3]{GP74} we have
\begin{proposition}
 \label{coincidence_GP}
  The elements $e_\alpha$, (resp. $f_{{\alpha}0}$) and
  $\tilde{e}_\alpha$ (resp. $\tilde{f}_{{\alpha}0}$)
  coincide in $D^4/\theta$.
\end{proposition}

Recall definitions of $\tilde{e}_\alpha$ and
$\tilde{f}_{{\alpha}0}$ from \cite{GP74}. \vspace{2mm}

\underline{The definition of $\tilde{e}_\alpha$,
\cite[p.6]{GP74}}.

\begin{equation}
  \label{GP_rec_def}
  \begin{split}
   & \tilde{e}_{i_n{i}_{n-1}\dots{i}_2{i}_1} =
     \tilde{e}_{i_n}\sum\limits_{\beta \in \Gamma_{e}(\alpha)}\tilde{e}_\beta,
  \end{split}
\end{equation}
where
\begin{equation}
  \begin{split}
   \Gamma_{e}(\alpha) = \{ \beta = (k_{n-1},\dots,k_2,k_1) \mid~
     & k_{n-1} \notin \{i_n,{i}_{n-1}\},\dots, k_1 \notin
     \{i_2,{i}_1\},\text{ and } \\
    &  k_1 \neq k_2, \dots, k_{n-2} \neq k_{n-1}  \}.
  \end{split}
\end{equation}

\underline{The definition of $\tilde{f}_{{\alpha}0}$,
\cite[p.53]{GP74}}.

\begin{equation}
  \begin{split}
   & \tilde{f}_{i_n{i}_{n-1}\dots{i}_2{i}_1{0}} =
     \tilde{e}_{i_n}\sum\limits_{\beta \in \Gamma_{f}(\alpha)}\tilde{e}_\beta,
  \end{split}
\end{equation}
where
\begin{equation}
  \begin{split}
   \Gamma_{f}(\alpha) = \{ \beta = (k_{n},\dots,k_2,k_1) \mid~
     & k_n \notin \{i_n,{i}_{n-1}\},\dots, k_2 \notin
     \{i_2,{i}_1\}, k_1 \notin \{i_1\} \text{ and } \\
    &  k_1 \neq k_2, \dots, k_{n-2} \neq k_{n-1}  \}.
  \end{split}
\end{equation}

\begin{proposition}[The elements $\tilde{e}_\alpha$]
  \label{coincidence_E}
 Consider elements $\tilde{e}_\alpha$ for $\alpha = 21, 121, 321,
2341$ (see Section \ref{examples_D4}). The relation
$$
   e_\alpha = \tilde{e}_\alpha
$$
takes place in $D^4$.
\end{proposition}
\PerfProof \underline{For $n=2$}: $\alpha = 21$. We have
\begin{equation}
  \tilde{e}_{21} = e_2\sum\limits_{j \neq 1,2}{e_j} =
    e_2(e_3 + e_4).
\end{equation}
According to Section \ref{examples_D4}, we see that $e_{21} =
\tilde{e}_{21}$.
 \vspace{2mm}

 \underline{For $n=3$}: 1) $\alpha = 121$,
\begin{equation}
  \begin{split}
  & \Gamma_{e}(\alpha) = \{ (k_2 k_1) \mid k_2 \in \{3,4 \}, k_1 \in
  \{3, 4\}, k_1 \neq k_2 \}, \\
  & \tilde{e}_{121} = e_1\sum\limits_{\beta \in \Gamma_{e}(\alpha)}e_\beta =
     e_1(e_{34} + e_{43}) = e_1(e_3(e_1 + e_2) + e_4(e_1 + e_2)) = e_1a^{34}_2.
  \end{split}
\end{equation}
By Section \ref{examples_D4} we have $e_{121} = \tilde{e}_{121}$.
\vspace{2mm}

 2) $\alpha = 321 = 341$. We have
\begin{equation*}
  \begin{split}
   \Gamma_{e}(\alpha) = & \{ (k_2 k_1) \mid k_2 \in \{3,2 \}, k_1 \in
  \{2, 1\}, k_1 \neq k_2 \} = \{ 14, 13, 43\}. \\
  \tilde{e}_{321} = & e_1\sum\limits_{\beta \in \Gamma_{e}(\alpha)}e_\beta =
     e_3(e_{14} + e_{13} + e_{43}) = \\
     & e_3(e_1(e_2 + e_4) + e_1(e_2 + e_3) + e_4(e_1 + e_2)) = \\
  &  e_3((e_1 + e_2)(e_2 + e_4)(e_1 + e_4) + e_1(e_2 + e_3)) = \\
  &  e_3((e_1 + e_2)(e_1 + e_4)(e_2 + e_4 + e_1(e_2 + e_3)) = \\
  &  e_3((e_1 + e_2)(e_1 + e_4)(e_4 + (e_1 + e_2)(e_2 + e_3)) = \\
  &  e_3((e_1 + e_2)(e_1 + e_4)(e_4(e_2 + e_3) + (e_1 + e_2)).
  \end{split}
\end{equation*}
Since $e_1 + e_2 \subseteq e_4(e_2 + e_3) + e_1 + e_2$, we have
\begin{equation}
   \tilde{e}_{321} = e_3(e_{14} + e_{13} + e_{43}) =
     e_3(e_1 + e_2)(e_1 + e_4).
\end{equation}
Since $\tilde{e}_{321}$ is symmetric with respect to transposition
$2 \leftrightarrow 4$, we have
\begin{equation}
   \tilde{e}_{321} = \tilde{e}_{341} =
   e_3(e_{14} + e_{13} + e_{43}) =
   e_3(e_{12} + e_{13} + e_{23}) =
     e_3(e_1 + e_2)(e_1 + e_4).
\end{equation}
By Section \ref{examples_D4} we have $e_{321} = \tilde{e}_{321}$.
 \vspace{2mm}

 \underline{For $n=4$}:  $\alpha = 2341 = 2321 = 2141$. We have
\begin{equation}
  \begin{split}
   \Gamma_{e}(\alpha) = & \{ (k_3 k_2 k_1) \mid k_3 \in \{1,4 \}, k_2 \in
  \{1, 2\}, k_1 \in \{ 2,3 \},\quad k_1 \neq k_2 , k_2 \neq k_3 \} = \\
  & \{ (123), (412), (413) = (423) \},
  \end{split}
\end{equation}
and
\begin{equation}
  \begin{split}
   \tilde{e}_{2341} = & e_2\sum\limits_{\beta \in \Gamma_{e}(\alpha)}e_\beta =
     e_2(e_{123} + e_{412} + e_{413}) = \\
   & e_2(e_1(e_2 + e_3)(e_4 + e_3) + e_4(e_1 + e_2)(e_3 + e_2) +
         e_4(e_1 + e_3)(e_2 + e_3)) = \\
   & e_2(e_4 + e_3)(e_1(e_2 + e_3) + e_4(e_1 + e_2)(e_3 + e_2) +
         e_4(e_1 + e_3)(e_3 + e_2)) =  \\
   & e_2(e_4 + e_3)(e_1(e_2 + e_3) + e_4 +
         e_4(e_1 + e_3)(e_1 + e_2)(e_3 + e_2)) =  \\
   & e_2(e_4 + e_3)(e_1(e_2 + e_3) + e_4).
  \end{split}
\end{equation}
By Section \ref{examples_D4} we have $e_{2341} =
\tilde{e}_{2341}$.
 The proposition is proved.
\qedsymbol \vspace{2mm}

\begin{proposition}[The elements $\tilde{f}_{{\alpha}0}$]
  \label{coincidence_F}
Consider elements $\tilde{f}_{{\alpha}0}$ for $\alpha = 21, 121,
321$ (see Section \ref{examples_D4}). The following relation
$$
   f_{{\alpha}0} = \tilde{f}_{{\alpha}0}
$$
takes place in $D^4$.
\end{proposition}
\PerfProof
 \underline{For $n=2$}: $\alpha = 21$. We have
\begin{equation}
 \begin{split}
  & \Gamma_{f}(\alpha) = \{ (k_2 k_1) \mid k_2 \in \{3,4 \}, k_1 \in
  \{2, 3, 4\}, k_1 \neq k_2 \}, \\
   \tilde{f}_{210} = & e_2\sum\limits_{\beta \in \Gamma_{f}(\alpha)}{e_\beta} =
    e_2(e_{32} + e_{34} + e_{42} + e_{43}) = \\
  & e_2(e_3(e_1 + e_4) + e_3(e_2 + e_1) +
    e_4(e_1 + e_3) + e_4(e_1 + e_2)) = \\
  & e_2((e_3 + e_4)(e_1 + e_4)(e_1 + e_3) + e_3(e_1 + e_2) + e_4(e_1 + e_2)) = \\
  & e_2(e_3 + e_4)((e_1 + e_4)(e_1 + e_3) + e_3(e_1 + e_2) + e_4(e_1 + e_2)) = \\
  & e_2(e_3 + e_4)(e_1 + e_4(e_1 + e_3) + e_3(e_1 + e_2) + e_4(e_1 + e_2)) = \\
  & e_2(e_3 + e_4)(e_4(e_1 + e_3) + (e_1 + e_3)(e_1 + e_2) + e_4(e_1 + e_2)) = \\
  & e_2(e_3 + e_4)(e_4(e_1 + e_3)(e_1 + e_2) + (e_1 + e_3)(e_1 + e_2) + e_4) = \\
  & e_2(e_3 + e_4)(e_4 + (e_1 + e_3)(e_1 + e_2)) = \\
  & e_2(e_3 + e_4)(e_4 + e_1 + e_3(e_1 + e_2)).
 \end{split}
\end{equation}
By Section \ref{examples_D4} we have
 $f_{210} = \tilde{f}_{210}$.
  \vspace{2mm}

\underline{For $n=3$}: 1) $\alpha = 121$.
 For this case, we have
\begin{equation}
   \Gamma_{f}(\alpha) = \{ (k_3 k_2 k_1) \mid k_3 \in \{3,4 \}, k_2 \in
  \{3, 4\}, k_1 \in \{2,3,4\}, k_1 \neq k_2, k_2 \neq k_3 \}, \\
\end{equation}
and
\begin{equation}
 \begin{split}
  \tilde{f}_{1210} = & e_1\sum\limits_{\beta \in \Gamma_{f}(\alpha)}{e_\beta} =
    e_1(e_{342} + e_{343} + e_{432} + e_{434}) = \\
    e_1[& e_3(e_1 + e_2)(e_4 + e_2) + e_3(e_1(e_4 + e_3) + e_2(e_4 + e_3))+ \\
       & e_4(e_1 + e_2)(e_3 + e_2) + e_4(e_1(e_4 + e_3) + e_2(e_4 + e_3))]= \\
    e_1[& e_3(e_1 + e_2)(e_4 + e_2 + e_3(e_1 + e_2(e_4 + e_3))) + \\
        & e_4(e_1 + e_2)(e_3 + e_2 + e_4(e_1 + e_2(e_4 + e_3)))] = \\
    e_1[& e_3(e_1 + e_2)(e_4 + e_2 + e_1(e_3 + e_2(e_4 + e_3))) + \\
        & e_4(e_1 + e_2)(e_3 + e_2 + e_1(e_4 + e_2(e_4 + e_3)))] = \\
    e_1[& (e_4 + e_2 + e_1(e_3 + e_2)(e_4 + e_3))]  \\
       [& (e_3 + e_2 + e_1(e_4 + e_2)(e_4 + e_3))]  \\
       [& e_3(e_1 + e_2) + e_4(e_1 + e_2)] = \\
    e_1[& (e_4 + e_2)(e_4 + e_3) + e_1(e_3 + e_2))]  \\
       [& (e_3 + e_2 + e_1(e_4 + e_2)(e_4 + e_3))]  \\
       [& e_3(e_1 + e_2) + e_4(e_1 + e_2)] = \\
    e_1[& (e_4 + e_2)(e_4 + e_3) + e_1(e_3 + e_2))]  \\
       [& e_1(e_3 + e_2) + (e_4 + e_2)(e_4 + e_3))]  \\
       [& e_3(e_1 + e_2) + e_4(e_1 + e_2)].
 \end{split}
\end{equation}
Since the first two intersection polynomials in the last
expression of $\tilde{f}_{1210}$ coincide with
$$
   (e_4 + e_2)(e_4 + e_3) + e_1(e_3 + e_2)),
$$
we have
\begin{equation}
   \tilde{f}_{1210} =
   e_1((e_4 + e_2)(e_4 + e_3) + e_1(e_3 + e_2)))(e_3 + e_4(e_1 + e_2)).
\end{equation}
By Section \ref{examples_D4} we have
 $f_{1210} = \tilde{f}_{1210}$.
  \vspace{2mm}

2) $\alpha = 321$. Here we have
\begin{equation}
 \begin{split}
  & \Gamma_{f}(\alpha) = \{ (k_3 k_2 k_1) \mid k_3 \in \{1,4 \}, k_2 \in
    \{3, 4\}, k_1 \in \{2,3,4\}, k_1 \neq k_2, k_2 \neq k_3 \} =  \\
  & \{ (132) = (142), (134), (143, (432), (434) \}
 \end{split}
\end{equation}
and
\begin{equation}
 \label{f_3210_transf_1}
 \begin{split}
   \tilde{f}_{3210} = & e_1\sum\limits_{\beta \in \Gamma_{f}(\alpha)}{e_\beta} =
    e_3(e_{132} + e_{134} + e_{143} + e_{432} + e_{434}) = \\
    e_3[& e_1(e_4 + e_2)(e_3 + e_2) + e_1(e_4 + e_2)(e_3 + e_4) + \\
        & e_1(e_3 + e_2)(e_3 + e_4) + e_4(e_3 + e_2)(e_1 + e_2) + \\
        & e_4(e_1 + e_2(e_3 + e_4))].
 \end{split}
\end{equation}
Since
\begin{equation}
 \begin{split}
  & e_1(e_4 + e_2)(e_3 + e_4) + e_4(e_1 + e_2(e_3 + e_4)) = \\
  & (e_3 + e_4)(e_4 + e_2)(e_1 + e_4(e_1 + e_2(e_3 + e_4)) = \\
  & (e_3 + e_4)(e_4 + e_2)(e_1 + e_4)(e_1 + e_2(e_3 + e_4)),
 \end{split}
\end{equation}
by (\ref{f_3210_transf_1}) we have
\begin{equation}
 \label{f_3210_transf_2}
 \begin{split}
   \tilde{f}_{3210} = & \\
    e_3[& e_1(e_4 + e_2)(e_3 + e_2) + \\
        & (e_3 + e_4)(e_4 + e_2)(e_1 + e_4)(e_1 + e_2(e_3 + e_4)) + \\
        &  e_1(e_3 + e_2)(e_3 + e_4) +
           e_4(e_3 + e_2)(e_1 + e_2)] = \\
    e_3(&e_1 + e_2)(e_1 + e_4)[e_1(e_4 + e_2)(e_3 + e_2) + \\
        & (e_3 + e_4)(e_4 + e_2)(e_1 + e_2(e_3 + e_4)) + \\
        & e_1(e_3 + e_2)(e_3 + e_4) +
          e_4(e_3 + e_2)] = \\
    e_3(&e_1 + e_2)(e_1 + e_4)[e_1(e_4 + e_2)(e_3 + e_2)(e_3 + e_4) + \\
        & (e_3 + e_4)(e_4 + e_2)(e_1 + e_2(e_3 + e_4)) + \\
        & e_1(e_3 + e_2) +
          e_4(e_3 + e_2)].
 \end{split}
\end{equation}
Since
$$
  e_1(e_4 + e_2)(e_3 + e_2)(e_3 + e_4) \subseteq e_1(e_3 + e_2),
$$
by (\ref{f_3210_transf_2}) we have
\begin{equation}
 \label{f_3210_transf_3}
 \begin{split}
   \tilde{f}_{3210} = & \\
    e_3(&e_1 + e_2)(e_1 + e_4)
        [(e_3 + e_4)(e_4 + e_2)(e_1 + e_2(e_3 + e_4)) + \\
        & e_1(e_3 + e_2) + e_4(e_3 + e_2)] = \\
    e_3(&e_1 + e_2)(e_1 + e_4)
        [e_1(e_3 + e_4)(e_4 + e_2) + e_2(e_3 + e_4)) + \\
        & e_1(e_3 + e_2) + e_4(e_3 + e_2)] = \\
    e_3(&e_1 + e_2)(e_1 + e_4)
        [e_1(e_3 + e_4)(e_4 + e_2)(e_3 + e_2) + e_2(e_3 + e_4)) + \\
        & e_1 + e_4(e_3 + e_2)].
 \end{split}
\end{equation}
Since
$$
   e_1(e_3 + e_4)(e_4 + e_2)(e_3 + e_2) \subseteq e_1,
$$
by (\ref{f_3210_transf_3}) we have
\begin{equation}
 \label{f_3210_transf_4}
 \begin{split}
   \tilde{f}_{3210} = & \\
    e_3(&e_1 + e_2)(e_1 + e_4)
        (e_2(e_3 + e_4)) + e_1 + e_4(e_3 + e_2)) = \\
    e_3(&e_1 + e_2)(e_1 + e_4)
        (e_1 + (e_2 + e_4)(e_3 + e_2)(e_3 + e_4)).
 \end{split}
\end{equation}

By Section \ref{examples_D4} we have
 $f_{3210} = \tilde{f}_{3210}$.
 The proposition is proved.
\qedsymbol \vspace{2mm}

\begin{conjecture}
 \label{conj_4}{\rm
  For every admissible sequence $\alpha$,
  the elements $e_\alpha$ (resp. $f_{{\alpha}0}$) and
  $\tilde{e}_\alpha$ (resp. $\tilde{f}_{{\alpha}0}$) coincide
  in $D^4$
  (see Proposition \ref{coincidence_GP})}.
\end{conjecture}

In Propositions \ref{coincidence_E} and \ref{coincidence_F}, this
conjecture was proven for small values of lengths of the
admissible sequence $\alpha$.

%% file: 4endom_D4.tex
 \section{\sc\bf Admissible elements in $D^4$ and Herrmann's
polynomials}
   \label{sect_Herrmann}

  In this section we consider Herrmann's endomorphisms $\gamma_{ik}$
  ($i,k = 1,2,3$) and polynomials $s_n$, $t_n$, $p_{i,n}$
  ($i = 1,2,3,4$) being perfect elements, \cite{H82}.
  Endomorphisms $\gamma_{ik}$ are important in
  construction of the perfect elements $s_n$, $t_n$, $p_{i,n}$,
  the perfect elements constitute $16$-element Boolean cube,
  coinciding modulo linear equivalence
  with the Gelfand-Ponomarev boolean cube $B^+(n)$, see
   Theorem \ref{th_Herrmann_82} due to Herrmann, \cite{H82}.

   We show that endomorphisms $\gamma_{ik}$ are also closely
 connected with admissible elements.

 First of all, the endomorphism $\gamma_{ik}$ acts on
 the admissible element $e_{\alpha{k}}$ such that
 \begin{equation*}
    \gamma_{ik}(e_{\alpha{k}}) = e_{\alpha{ki}},
 \end{equation*}
 see Theorem \ref{relat_rho_ik}. Thus, endomorphism $\gamma_{ik}$
 acts also on the admissible sequence. We see some similarity
 between the action of endomorphism $\gamma_{ik}$ and the action of
 elementary map of Gelfand-Ponomarev $\phi_i$.
 The endomorphism $\gamma_{ik}$ and the elementary map $\phi_i$ act in a sense
 in opposite directions, namely the endomorphism $\gamma_{ik}$ adds the index
 to the \underline{start} of the admissible sequence,
 and the elementary map $\phi_i$ adds
 the index to the \underline{end}\footnote{Recall, that in the admissible
 sequence $i_n{i}_{n-1}\dots{i}_2{i}_1$ the index $i_1$ is the start and $i_n$
 is the end.}.

 Further, the endomorphisms $\gamma_{ik}$ commute and we can consider
 a sequence of these endomorphism. Admissible elements
 $e_1{a}^{34}_t{a}^{24}_s{a}^{32}_r$
 from Table \ref{tab_gen_form_D4_e}
 are obtained by means of Herrmann's endomorphisms as follows:
 \begin{equation*}
      \gamma^t_{12}\gamma^s_{13}\gamma^r_{14}(e_1) =
             e_1{a}^{34}_t{a}^{24}_s{a}^{32}_r,
 \end{equation*}
   see Theorem \ref{theor_adm_Herrmann}.

 At last, we will see how Herrmann's
 polynomials $s_n$, $t_n$, $p_{i,n}$ are expressed by means
 of the cumulative elements $e_i(n)$, $f_i(n)$:
 \begin{equation*}
   \begin{split}
     & s_n = \sum\limits_{i=1,2,3,4}e_i(n), \\
     & t_n = f_0(n+1), \\
     & p_{i,n} = e_i(n) + f_0(n+1), \\
   \end{split}
 \end{equation*}
  see Theorem \ref{prop_endom_R} and Table \ref{tab_perfect_in D4}.

 \subsection{An unified formula of admissible elements}

  It turned out, that the admissible elements $e_\alpha$ and $e_{\alpha0}$
  described by Table \ref{table_adm_elem_D4} can be written by an unified formula.
  This unified formula represented in Table \ref{tab_gen_form_D4_e} and
  Proposition \ref{prop_gen_form_D4}.

  A main difference between Table \ref{table_adm_elem_D4} and
  Table \ref{tab_gen_form_D4_e} is in following:
  Table \ref{table_adm_elem_D4} describes admissible elements
  with admissible sequence $\alpha$
  \underline{starting} at generator $e_1$, and
  Table \ref{tab_gen_form_D4_e} describes
  admissible elements with admissible sequence $\alpha$
  \underline{ending} at generator $e_1$.
  Recall, that type $Fij$ (resp. $Gij$, $H11$) denotes the
  admissible sequence starting at $j$ and ending at $i$,
  see Remark \ref{def_Fij}.

  As we will see in the next sections,
  the unified formula represented in this section in
  Table \ref{tab_gen_form_D4_e} and in Proposition
  \ref{prop_gen_form_D4} is a basis of the construction
  of Herrmann's polynomials,
  \cite{H82}, \cite{H84}.

 \begin{table} [h]
   \renewcommand{\arraystretch}{1.5}
  \begin{tabular} {||c|c|c|c||}
  \hline \hline
    & Admissible  & Admissible  & Signature of \\
    & polynomial $e_\alpha$  & polynomial $f_{\alpha0}$
           & polynomial $e_\alpha$ \\
   \hline \hline
    $F12$ & $e_1{a}^{32}_{2r}a^{42}_{2s}a^{34}_{2t-1}$
        & $e_\alpha(e_1a^{34}_{2t} + a^{42}_{2s+1}{a}^{32}_{2r-1})$
        & 0 \quad 0 \quad 1  \\
  \hline
    $F13$ & $e_1{a}^{32}_{2r}a^{42}_{2s-1}a^{34}_{2t}$
        & $e_\alpha(e_1a^{34}_{2t+1} + a^{42}_{2s}{a}^{32}_{2r-1})$
        & 0 \quad 1 \quad 0  \\
  \hline
    $F14$ & $e_1{a}^{32}_{2r-1}a^{42}_{2s}a^{34}_{2t}$
        & $e_\alpha(e_1a^{34}_{2t+1} + a^{42}_{2s+1}{a}^{32}_{2r-2})$
        & 1 \quad 0 \quad 0  \\
  \hline \hline
    $G11$ & $e_1{a}^{32}_{2r}a^{42}_{2s}a^{34}_{2t}$
        & $e_\alpha(e_1a^{34}_{2t+1} + a^{42}_{2s+1}{a}^{32}_{2r-1})$
        & 0 \quad 0 \quad 0  \\
  \hline
    $G12$ & $e_1{a}^{32}_{2r+1}a^{42}_{2s-1}a^{34}_{2t}$
        & $e_\alpha(e_1a^{34}_{2t+1} + a^{42}_{2s}{a}^{32}_{2r})$
        & 1 \quad 1 \quad 0  \\
  \hline
    $G13$ & $e_1{a}^{32}_{2r-1}a^{42}_{2s}a^{34}_{2t+1}$
        & $e_\alpha(e_1a^{34}_{2t} + a^{42}_{2s+1}{a}^{32}_{2r-2})$
        & 1 \quad 0 \quad 1  \\
  \hline
    $G14$ & $e_1{a}^{32}_{2r}a^{42}_{2s-1}a^{34}_{2t+1}$
        & $e_\alpha(e_1a^{34}_{2t+2} + a^{42}_{2s}{a}^{32}_{2r-1})$
        & 0 \quad 1 \quad 1  \\
  \hline \hline
    $H11$ & $e_1{a}^{32}_{2r-1}a^{42}_{2s+1}a^{34}_{2t-1}$
        & $e_\alpha(e_1a^{34}_{2t} + a^{42}_{2s+2}{a}^{32}_{2r-2})$
        & 1 \quad 1 \quad 1  \\
  \hline \hline
  \end{tabular}
   \vspace{2mm}
   \caption{\hspace{3mm}Signatures of admissible elements
      $e_\alpha$ and $f_{\alpha0}$ ending at $1$}
   \label{tab_gen_form_D4_e}
 \end{table}

 \begin{proposition}
   \label{prop_gen_form_D4}
  {\rm(a)} The set of all admissible polynomials $e_\alpha$
    (described by Table \ref{table_adm_elem_D4})
   with admissible sequences $\alpha$ ending at $1$ coincides
   with the set of polynomials
   \begin{equation}
    \begin{split}
     \label{gen_from_D4}
       & e_1{a}_r^{32}{a}_s^{24}{a}_t^{34}, \quad \text{ where } r,s,t = 0,1,2,\dots, \\
     \end{split}
   \end{equation}
   see Table \ref{tab_gen_form_D4_e}.

   Similarly, the set of all admissible polynomials $e_\alpha$
   with admissible sequences $\alpha$ ending at $i \in \{2,3,4\}$
   coincides with the set of polynomials
   \begin{equation}
    \begin{split}
     \label{gen_from_D4_234}
        e_i{a}_r^{jl}{a}_s^{lk}{a}_t^{kj},
           \quad & \text{ where } r,s,t = 0,1,2,\dots, \\
       &  \text{ and } \{i,j,k,l\} = \{1,2,3,4\}.
     \end{split}
   \end{equation}

  {\rm(b)} The set of all admissible polynomials $f_{\alpha0}$,
    (described by Table \ref{table_adm_elem_D4})
   with admissible sequences $\alpha$ ending at $1$ coincides
   with the set of polynomials
   \begin{equation}
    \begin{split}
     \label{gen_from_D4_f}
       e_\alpha(e_1{a}_{t+1}^{34}+ {a}_{s+1}^{42}a_{r-1}^{32}) & =
         e_1{a}_t^{34}{a}_s^{42}a_r^{32}
           (e_1{a}_{t+1}^{34}+ {a}_{s+1}^{42}a_{r-1}^{32}), \\
        \text{ where }& r,s,t = 0,1,2,\dots,
     \end{split}
   \end{equation}
   see Table \ref{tab_gen_form_D4_e}.

   Similarly, the set of all admissible polynomials $f_{\alpha0}$
   with admissible sequences $\alpha$ ending at $i \in \{2,3,4\}$
   coincides with the set of polynomials
   \begin{equation}
    \begin{split}
     \label{gen_from_D4_234_f}
       e_\alpha(e_i{a}_t^{jk}+ {a}_s^{kl}a_r^{jl}) = &
         e_i{a}_t^{jk}{a}_s^{kl}a_r^{lj}
           (e_i{a}_{t+1}^{jk}+ {a}_{s+1}^{kl}a_{r-1}^{jl}), \\
        & \text{ where } r,s,t = 0,1,2,\dots,
         \text{ and } \{i,j,k,l\} = \{1,2,3,4\}.
     \end{split}
   \end{equation}
 \end{proposition}

  \PerfProof  For any atomic element $a^{pq}_i$, set the range $1$ if $i$ is
  odd and $0$ if $i$ is even. For the admissible polynomials
  $e_\alpha$ from Table \ref{tab_gen_form_D4_e},
  the set of the corresponding ranges we call the {\it signature}.
  For example, for $e_\alpha$ of the type $F12$, the
  signature is $(0, 0, 1)$, see Table \ref{tab_gen_form_D4_e}.
  There are $8$ possible signatures.
  Since, all possible signatures appeared in
  Table \ref{tab_gen_form_D4_e},
  we can get all possible combination of indices $(r, s, t)$
  in (\ref{gen_from_D4}), (\ref{gen_from_D4_f}). \qedsymbol

  \begin{remark}
   \label{rem_indices_D4}
  {\rm  For the admissible polynomials $e_\alpha$ the sum of low
  indices of atomic elements $a^{kl}_r$ is by $1$ less the length
  of the admissible sequence $\alpha$,
  for the admissible polynomials $f_{\alpha0}$ the sum of low
  indices of atomic elements $a^{kl}_r$ (contained in the parantheses)
  is equal to the length of the admissible sequence
  $\alpha$, see Table \ref{table_adm_elem_D4}.
  }
  \end{remark}

 \subsection{Inverse cumulative elements in $D^4$}
   \label{sect_inverse_cum}
   In addition to cumulative elements, we introduce now {\it inverse
   cumulative elements}. A difference between cumulative elements
   and inverse cumulative elements is in following: the cumulative elements
   accumulate all admissible elements of the given length
   \underline{starting} at some generator $e_i$, and
   inverse cumulative elements accumulate all admissible elements of the given length
   \underline{ending} at some generator $e_i$.

   The cumulative elements $e_i(n)$, where $i = 1,2,3,4$, are defined in
   (\ref{cumul_polyn_D4}). Proposition \ref{prop_gen_form_D4}
   and Remark \ref{rem_indices_D4} motivate the following definition
   of inverse cumulative polynomials $e_i^\vee(n)$, where $i = 1,2,3,4$,
   and $f_0^\vee(n)$ as follows:

   \begin{equation}
    \label{cumul_polyn_e}
     \begin{split}
         & e_1^\vee(n) = \sum\limits_{r + s + t = n-1}
                      {e_1{a}_r^{32}{a}_s^{24}{a}_t^{34}},
             \qquad \qquad
           e_2^\vee(n) = \sum\limits_{r + s + t = n-1}
                      {e_2{a}_r^{31}{a}_s^{14}{a}_t^{34}}, \\
         & \\
         & e_3^\vee(n) = \sum\limits_{r + s + t = n-1}
                      {e_3{a}_r^{12}{a}_s^{24}{a}_t^{14}},
             \qquad \qquad
           e_4^\vee(n) = \sum\limits_{r + s + t = n-1}
                      {e_4{a}_r^{12}{a}_s^{23}{a}_t^{13}},
     \end{split}
   \end{equation}
 and
   \begin{equation}
    \label{cumul_polyn_f}
     \begin{split}
         f_0^\vee(n) = & \sum f_{1{i}_{n-1}\dots{i}_2{0}} +
                    \sum f_{2{i}_{n-1}\dots{i}_2{0}} +
                    \sum f_{3{i}_{n-1}\dots{i}_2{0}} +
                    \sum f_{4{i}_{n-1}\dots{i}_2{0}} = \\
         & \\
         & \sum\limits_{r + s + t = n}
                      e_\alpha({e_1{a}_t^{34} + {a}_s^{24}a_r^{32}}) +
           \sum\limits_{r + s + t = n}
                      e_\alpha({e_2{a}_t^{34} + {a}_s^{14}a_r^{31}})+ \\
         & \\
         & \sum\limits_{r + s + t = n}
                      e_\alpha({e_3{a}_t^{24} + {a}_s^{14}a_r^{21}}) +
           \sum\limits_{r + s + t = n}
                      e_\alpha({e_4{a}_t^{23} + {a}_s^{13}a_r^{21}}).
     \end{split}
   \end{equation}

 \begin{proposition}
   \label{sum_of_all}
      {\rm(a)} The sum of all cumulative elements $e_i(n)$
      of the given length $n$ and the sum of all
      inverse cumulative elements $e_i^{\vee}(n)$
      coincide, i.e.,
     \begin{equation}
          e_1(n) + e_2(n) + e_3(n) + e_4(n) =
          e_1^{\vee}(n) + e_2^{\vee}(n) + e_3^{\vee}(n) +
          e_4^{\vee}(n).
     \end{equation}

     {\rm(b)} The cumulative element $f_0(n)$ coincides with inverse
     cumulative element $f_0^{\vee}(n)$:
     \begin{equation}
          f_0(n) = f_0^{\vee}(n).
     \end{equation}

     {\rm(c)} The number of elements in every sum (\ref{cumul_polyn_e}) is
     $\displaystyle\frac{1}{2}(n+1)(n+2)$\footnote{Compare with Remark \ref{slices}.}.
 \end{proposition}
   \PerfProof (a), (b) are true since sums from the both sides
   consist of all admissible elements of the
   given length $n$.

\begin{figure}[h]
 \includegraphics{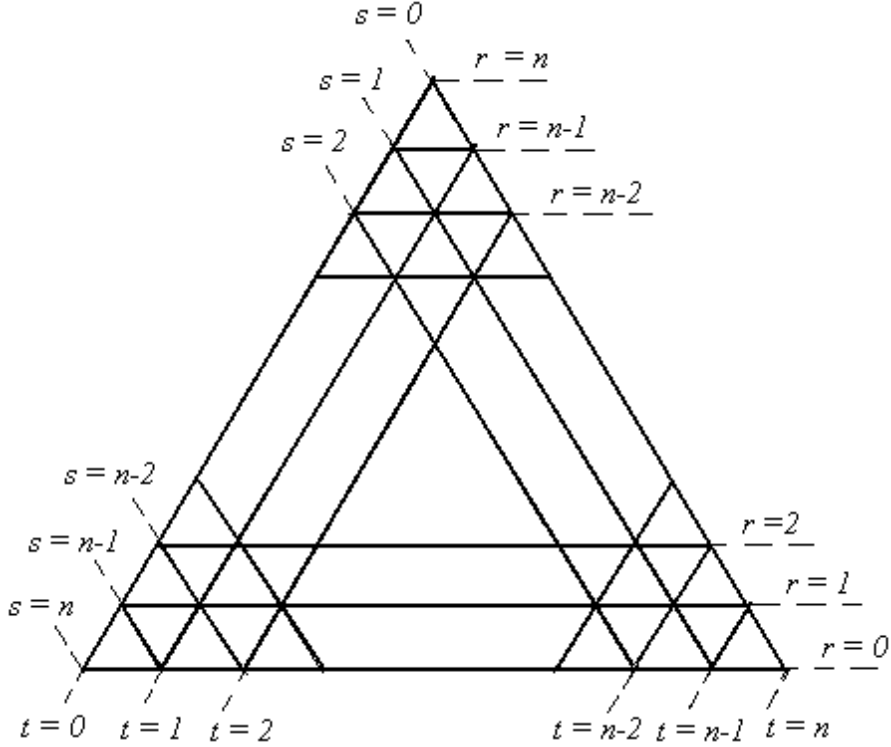}
 \caption{\hspace{3mm}The triangle with integer barycentric
 coordinates}
 \label{barycentric_coord}
 \end{figure}

   (c) We just need to found the number of solutions of the equation
   \begin{equation}
          r + s + t = n.
   \end{equation}
   These solutions are points with integer
   {\it barycentric coordinates}\footnote{For details concerning barycentric coordinates,
   see, e.g., H.~S.~M.~Coxeter's book, \cite[Section 13.7]{Cox89}, or
   A.~Bogomolny's site \cite{Bm96}.} in the triangle depicted on
   Fig. \ref{barycentric_coord}. Let $(r, s, t)$ be coordinates of
   any point of this triangle. Any move along
   one of edges does not change the coordinate sum
   $r + s + t$, and this sum is equal to $n$.
   \qedsymbol
   ~\vspace{3mm} \\

  \subsection{Herrmann's endomorphisms and admissible elements}
    \label{sect_Herrmann_introduced}
  Herrmann introduced in \cite[p.361, p.367]{H82}, \cite[p.229]{H84}
  polynomials $q_{ij}$ and associated endomorphisms
  $\gamma_{ij}$ of $D^4$ playing the central
  role in his study of the modular lattice
  $D^4$, in particular, in his construction of perfect polynomials.

  For $\{i,j,k,l\} = \{1,2,3,4\}$, define
  \begin{equation}
    \label{poly_qij}
     q_{ij} = q_{ji} = q_{kl} = q_{lk} = (e_i + e_j)(e_k + e_l).
  \end{equation}
 In our denotations,
  \begin{equation}
    \label{poly_qij_our}
     q_{ij} = a_1^{ij}a_1^{kl}, \text{ where } \{i,j,k,l\} = \{1,2,3,4\}.
  \end{equation}
  The endomorphism $\gamma_{ij}$ of $D^4$ is denoted as follows:
 \begin{equation}
    \label{endom_rho_ij}
     1 \mapsto q_{ij}, \quad 0 \mapsto 0, \quad e_k \mapsto e_k{q_{ij}}.
  \end{equation}
 For every polynomial $f(e_1,e_2,e_3,e_4)$, we have
 \begin{equation*}
    \gamma_{ij}f(e_1,e_2,e_3,e_4) =
    f(e_1{q}_{ij},e_2{q}_{ij},e_3{q}_{ij},e_4{q}_{ij}).
 \end{equation*}
 Essentially, by (\ref{poly_qij}) among endomorphisms $\gamma_{ij}$,
 there are only $3$ different:
 \begin{equation*}
    \gamma_{12}, \quad  \gamma_{13}, \quad \gamma_{14}.
 \end{equation*}

 \begin{proposition}
    All Herrmann's endomorphisms $\gamma_{ij}$ commute:
    \begin{equation}
       \gamma_{1i}\gamma_{1j} = \gamma_{1j}\gamma_{1i}, \quad
          i, j \in \{2,3,4\}, \quad i \neq j.
    \end{equation}
 \end{proposition}

 \PerfProof It suffices to check the commutativity on generators. We
 will check that
    \begin{equation}
       \gamma_{13}(\gamma_{12}(e_i)) = \gamma_{12}(\gamma_{13}(e_i)), \quad
        \text{ where } i = 1,2,3,4.
    \end{equation}
 We have $\gamma_{12}(e_1) = e_1(e_3 + e_4)$, and
 \begin{equation*}
  \begin{split}
    \gamma_{13}(\gamma_{12}(e_1)) = & \gamma_{13}(e_1(e_3 + e_4)) =
     \gamma_{13}(e_1)(\gamma_{13}(e_3) + \gamma_{13}(e_4)) = \\
     & e_1(e_2 + e_4)(e_3(e_2 + e_4) + e_4(e_1 + e_3)) = \\
     & e_1(e_2 + e_4)(e_3 + e_4)(e_1 + e_3) = \\
     & e_1(e_2 + e_4)(e_3 + e_4),
  \end{split}
 \end{equation*}
 and
 \begin{equation}
  \label{rho_13}
     \gamma_{13}(\gamma_{12}(e_1)) = \gamma_{12}(\gamma_{13}(e_1)).
 \end{equation}
 From (\ref{rho_13}) we have
 \begin{equation*}
   \begin{split}
    & \gamma_{13}(\gamma_{12}(e_2)) = \gamma_{24}(\gamma_{21}(e_2)) =
     \gamma_{21}(\gamma_{24}(e_2)) = \gamma_{12}(\gamma_{13}(e_2)), \\
    & \gamma_{13}(\gamma_{12}(e_3)) = \gamma_{31}(\gamma_{34}(e_3)) =
     \gamma_{34}(\gamma_{31}(e_3)) = \gamma_{12}(\gamma_{13}(e_3)), \\
    & \gamma_{13}(\gamma_{12}(e_4)) = \gamma_{42}(\gamma_{43}(e_4)) =
     \gamma_{43}(\gamma_{42}(e_4)) = \gamma_{12}(\gamma_{13}(e_4)).
     \qed
   \end{split}
 \end{equation*}

  \subsubsection{More relations on the admissible sequences}

  Further, we want to found a connection between Herrmann's
  endomorphisms $\gamma_{1i}$ and admissible sequences in $D^4$.
  As we will see in Proposition \ref{relat_rho_ik},
  endomorphisms $\gamma_{1i}$ add corresponding indices to the
  \underline{start} of the given admissible sequence, while elementary maps
  $\phi_i$ of Gelfand-Ponomarev add indices to the \underline{end} of
  the corresponding admissible sequence,
  see Table \ref{table_admissible_ExtD4}. To found this connection,
  we need more relations connecting admissible sequences in $D^4$,
  see Proposition \ref{relations_1_14} and the fundamental property
  of the admissible sequences (\ref{main_D4}).
 \begin{proposition}
   \label{more_relat_adm}
   The following relations hold
   \item[1)] $2(13)^s{1} = 2(42)^s{1}$,
   \item[2)] $2(13)^s(14)^r{1} = 2(42)^s(32)^r{1}$,
   \item[3)] $1(41)^r(31)^s(21)^t = (14)^r(13)^s(12)^t{1} =
             (12)^t(13)^s(14)^r{1}$,
   \item[4)] $1(41)^r(31)^s(21)^t = (12)^t(42)^s(32)^r{1}$,
   \item[5)] $3(41)^s2 = (32)^s12$,
   \item[6)] $3(41)^s(21)^t(31)^r = 3(41)^s(31)^t(21)^r =
                (32)^s(42)^r(12)^t1$,
   \item[7)] $2(42)^r(32)^s(12)^t1 = (21)^{t+1}(41)^s(31)^r$,
   \item[8)] $3(42)^s(12)^t(32)^r1 = (31)^s(41)^r(21)^{t+1}$,
   \item[9)] $1(42)^r(32)^s(12)^t1 = (13)^r(41)^{s+1}(21)^t$,
   \item[10)] $(23)^s(42)^r(12)^t1 = 2(41)^s(31)^{r-1}(21)^{t+1}$.
 \end{proposition}
  \PerfProof
  1) We have
 \begin{equation*}
  \begin{split}
     2(13)^s{1} = & 2(13)(13)(13)\dots(13)(13)1 =
       2(42)(13)(13)\dots(13)(13)1 = \\
     & 2(42)(42)(13)\dots(13)(13)1 =
       2(42)(42)(42)\dots(42)(42)1 = 2(42)^s{1}.
  \end{split}
 \end{equation*}
   2) By heading 1) we have
 \begin{equation*}
  \begin{split}
     2(13)^s(14)^r{1} = & 2(42)^s(32)^r{1} =
       2(42)^s(14)(14)(14)\dots(14)(14)1 = \\
     & 2(42)^s(32)(14)(14)\dots(14)(14)1 =
       2(42)^s(32)(32)(32)\dots(32)(32)1 = \\
     & 2(42)^s(32)^r1.
  \end{split}
 \end{equation*}
    3) First,
 \begin{equation*}
  \begin{split}
      1(41)^r = & 1(41)(41)\dots(41) = (14)(14)\dots(14)1 =
      (14)^r1, \text{ and } \\
      & 1(41)^r(31)^s(21)^t = (14)^r(13)^s(12)^t{1}.
      \end{split}
 \end{equation*}
  By heading 6) of Proposition \ref{main_D4}
 \begin{equation*}
  \begin{split}
      & 1(41)^r(31)^s(21)^t = 1(21)^t(31)^s(41)^r, \text{ and } \\
      & 1(41)^r(31)^s(21)^t = (12)^t(13)^s(14)^r{1}.
      \end{split}
 \end{equation*}
     4) By heading 3) and 2) we have
 \begin{equation*}
  \begin{split}
       1(41)^r(31)^s(21)^t = &
        (12)^t(13)^s(14)^r1 = (12)^{t-1}1[2(13)^s(14)^r{1}] = \\
      & (12)^{t-1}1[2(42)^s(32)^r{1}] = (12)^{t-1}12(42)^s(32)^r{1} = \\
      & (12)^{t}(42)^s(32)^r{1}.
      \end{split}
 \end{equation*}
     5) Here,
 \begin{equation*}
  \begin{split}
      3(41)^s2 = & 3(41)(41)\dots(41)(41)2 =
        3(23)(41)\dots(41)(41)2 = \\
      & 3(23)(23)\dots(23)212 = (32)^s12.
      \end{split}
 \end{equation*}
     6) By heading 5) of this proposition and
     by heading 6) of Proposition \ref{main_D4}we have
 \begin{equation*}
  \begin{split}
      3(41)^s(21)^t(31)^r = & (32)^s1(21)^t(31)^r = \\
      & (32)^s1(31)^r(21)^t = (32)^s(13)^r(12)^t1 = \\
      & (32)^s(13)(13)\dots(13)(13)(12)^t1 = \\
      & (32)^s(42)(13)\dots(13)(13)(12)^t1 = \\
      & (32)^s(42)^r(12)^t1.
      \end{split}
 \end{equation*}
     7) By heading 2)
 \begin{equation}
  \label{relat_7_adm}
  \begin{split}
      & 2(32)^s(42)^r1 = 2(14)^s(13)^r1 = \\
      & 2(13)^r(14)^s1 = 21(31)^r(41)^s. \\
      \end{split}
 \end{equation}
 By (\ref{relat_7_adm}) and by heading 6) of Proposition
 \ref{main_D4} we have
 \begin{equation*}
  \begin{split}
      & 2(42)^r(32)^s(12)^t1 = 2(12)^t(32)^s(42)^r1 = \\
      & 2(12)^t1(31)^r(41)^s = (21)^{t+1}(31)^r(41)^s = \\
      & (21)^{t+1}(41)^s(31)^r.
      \end{split}
 \end{equation*}
     8) By heading 5) and substitution $1 \leftrightarrow 2$
 \begin{equation}
  \label{relat_8_adm}
       3(42)^s1 = (31)^s21.
 \end{equation}
 By (\ref{relat_8_adm}) and by heading 1) of Proposition
 \ref{main_D4} we have
 \begin{equation*}
  \begin{split}
      3(42)^s(12)^t(32)^r1 = & (31)^s2(12)^t(32)^r1 =
       (31)^s(21)^t(23)^r21 = \\
      & (31)^s(23)^r(21)^{t+1} =
        (31)^s(23)(23)\dots(23)(23)(21)^{t+1} = \\
      & (31)^s(41)(23)\dots(23)(23)(21)^{t+1} =
        (31)^s(41)^r(21)^{t+1}. \\
      & \\
      \end{split}
 \end{equation*}
     9) By substitution $3 \rightarrow 1 \rightarrow 2 \rightarrow 3$
     in heading 5) we have
 \begin{equation}
  \label{relat_9_adm}
     1(42)^r3 = (13)^r23.
 \end{equation}
 In addition,
 \begin{equation}
  \label{relat_9_adm_1}
      \begin{split}
      3(23)^s21 = & 3(23)(23)\dots(23)(23)21 = 3(23)(23)\dots(23)(23)41 \\
       & 3(23)(23)\dots(23)(41)41 = 3(41)(41)\dots(41)(41)41 = 3(41)^{s+1}.
      \end{split}
 \end{equation}
   From (\ref{relat_9_adm}) and (\ref{relat_9_adm_1}) we get
 \begin{equation}
  \label{relat_9_adm_2}
      \begin{split}
       1(42)^r(32)^s(12)^t1 = & (13)^r2(32)^s(12)^t1 =
          (13)^r(23)^s2(12)^t1 = \\
       & (13)^r(23)^s(21)(21)^t = (13)^r(41)^{s+1}(21)^t.
      \end{split}
 \end{equation}
  10)
     First, we have
 \begin{equation}
  \label{relat_10_adm}
      \begin{split}
       (23)^s42 = & 2(32)(32)\dots(32)342 = 2(32)(32)\dots(32)412 =  \\
       & 2(41)(41)\dots(41)412 = 2(41)^s2,
      \end{split}
 \end{equation}
 and
 \begin{equation}
  \label{relat_10_adm_1}
      \begin{split}
       1(24)^{r-1}21 = & 1(24)(24)\dots(24)21 = 1(31)(24)\dots(24)21 =  \\
       & 1(31)(31)\dots(31)21 = 1(31)^{r-1}21.
      \end{split}
 \end{equation}
   From (\ref{relat_10_adm}) and (\ref{relat_10_adm_1}) we get
 \begin{equation}
  \label{relat_10_adm_2}
      \begin{split}
       (23)^s(42)^r(12)^t1 = & (23)^s4(24)^{r-1}2(12)^t1 = \\
       & 2(41)^s(24)^{r-1}2(12)^t1 = 2(41)^s(31)^{r-1}2(12)^t1 = \\
       & 2(41)^s(31)^{r-1}(21)^{t+1}.
      \end{split}
 \end{equation}
 The proposition is proved. \qedsymbol

 \subsubsection{How Herrmann's endomorphisms act on the admissible elements}

   Before understanding how Herrmann's endomorphisms act on the admissible
   elements, we should know how these endomorphisms work on the
   simple admissible elements $e_k{a}^{ij}$
   are {\it almost atomic} elements\footnote{Recall, that
   polynomials ${a}^{ij}_r$ are called {\it atomic}, see
   (\ref{def_atomic}).}.

 \begin{proposition}
   \label{relat_rho_1j}
    Endomorphisms $\gamma_{ij}$ act on admissible elements as follows:
  \begin{enumerate}
    \item $\gamma_{1j}(e_1{a}^{kl}_r) = e_1{a}^{kl}_{r+1}$,
          $\{1,j,k,l\} = \{1,2,3,4\}$, \vspace{3mm}
    \item $\gamma_{1k}(e_1{a}^{kl}_r) = \gamma_{jl}(e_1{a}^{kl}_r) =
          e_1{a}^{jl}_1{a}^{kl}_r$,
          $\{1,j,k,l\} = \{1,2,3,4\}$. \vspace{3mm}
   \end{enumerate}
 \end{proposition}

  \PerfProof
  Without lost of generality, we will show that
  \begin{equation}
    \label{prove_rho_1}
      \gamma_{12}(e_1{a}^{34}_r) = e_1{a}^{34}_{r+1},
  \end{equation}
  and
  \begin{equation}
    \label{prove_rho_1A}
       \gamma_{12}(e_1{a}^{23}_r) = e_1{a}^{34}_1{a}^{23}_{r}. \\
  \end{equation}

  1) Let us prove (\ref{prove_rho_1}). For r=0, we get
  \begin{equation}
    \label{base_rho12}
       \gamma_{12}(e_1) = e_1(e_3 + e_4) = e_1{a}^{34}_1.
  \end{equation}
  For r=1,
  \begin{equation*}
   \begin{split}
     \gamma_{12}(e_1{a}^{34}_1) = & \gamma_{12}(e_1(e_3 + e_4)) =
       \gamma_{12}(e_1)(\gamma_{12}(e_3) + \gamma_{12}(e_4)) = \\
     & e_1(e_3 + e_4)(e_3(e_1 + e_2) + e_4(e_1 + e_2)) =
       e_1(e_3(e_1 + e_2) + e_4(e_1 + e_2)) = \\
     & e_1(e_3 + e_4(e_1 + e_2)) = e_1a^{34}_2.
   \end{split}
  \end{equation*}
  By induction hypothesis and (\ref{poly_qij}) we have
  \begin{equation*}
   \begin{split}
      \gamma_{12}(e_1{a}^{34}_r) = & \gamma_{12}(e_1(e_3 + e_4{a}^{12}_{r-1})) =
       \gamma_{12}(e_1)(\gamma_{12}(e_3) + \gamma_{12}(e_4)\gamma_{34}({a}^{12}_{r-1})) = \\
     & e_1(e_3 + e_4)(e_3(e_1 + e_2) + e_4(e_1 + e_2){a}^{12}_{r})= \\
     & e_1(e_3 + e_4)(e_3 + e_4(e_1 + e_2){a}^{12}_{r}) = \\
     & e_1(e_3 + e_4)(e_3 + e_4{a}^{12}_{r}) =
       e_1{a}^{34}_{r+1}.
   \end{split}
  \end{equation*}

 2) Now, let us prove (\ref{prove_rho_1A}). For r=1,
  \begin{equation*}
   \begin{split}
     \gamma_{12}(e_1{a}^{23}_1) = & \gamma_{12}(e_1(e_2 + e_3)) =
       \gamma_{12}(e_1)(\gamma_{12}(e_2) + \gamma_{12}(e_3)) = \\
     & e_1(e_3 + e_4)(e_2(e_3 + e_4) + e_3(e_1 + e_2)) = \\
     & e_1(e_2(e_3 + e_4) + e_3(e_1 + e_2)) = \\
     & e_1(e_3 + e_4)(e_2 + e_3)(e_1 + e_2) =
       e_1(e_3 + e_4)(e_2 + e_3)) = \\
     & e_1{a}^{34}_1{a}^{23}_1.
   \end{split}
  \end{equation*}

 By induction hypothesis, and again by (\ref{poly_qij}) we have
  \begin{equation*}
   \begin{split}
     \gamma_{12}(e_1{a}^{23}_r) = & \gamma_{12}(e_1(e_2 + e_3{a}^{41}_{r-1})) =
       \gamma_{12}(e_1)(\gamma_{12}(e_2) + \gamma_{34}(e_3{a}^{41}_{r-1})) = \\
     & e_1(e_3 + e_4)(e_2(e_3 + e_4) + e_3{a}^{12}_1{a}^{41}_{r-1}) = \\
     & e_1(e_3 + e_4)(e_2 + e_3{a}^{12}_1{a}^{41}_{r-1}) = \\
     & e_1{a}^{34}_1(e_2 + e_3{a}^{41}_{r-1}) =
     e_1{a}^{34}_1{a}^{23}_{r}.
   \end{split}
  \end{equation*}
   The proposition is proved. \qedsymbol
 \begin{theorem}
   \label{relat_rho_ik}
    Endomorphism $\gamma_{ik}$ acts
     on the admissible element $e_{\alpha{k}}$ as follows:
       $$
         \gamma_{ik}(e_{i_{r}i_{r-1}\dots{i_2}{k}}) =
          e_{i_{r}i_{r-1}\dots{i_2}{ki}},
       $$
     or, in other words,
       \begin{equation}
         \label{lab_rho_ik}
          \fbox{$
          \gamma_{ik}(e_{\alpha{k}}) = e_{\alpha{ki}},$}
       \end{equation}
     where $i, k \in \{1,2,3,4\}$, $i \neq k$, and
       $\alpha, \alpha{k}$ are admissible sequences.
 \end{theorem}

  Without lost of generality, we will show that
  \begin{equation}
    \label{prove_rho_2}
       \qquad \gamma_{12}(e_{\alpha{2}}) = e_{\alpha{21}}.
  \end{equation}

   It suffices to prove (\ref{prove_rho_2}) for all cases $\alpha2$
  described by Table \ref{table_adm_elem_D4}
  with the permutation $1 \leftrightarrow 2$:
  $$
     F12, F32, F42, G22, G12, G32, G42, H22.
  $$

  \underline{Case $F12$}. Here,
  \begin{equation}
    \begin{split}
     & \alpha2 = (12)^t(42)^r(32)^s, \quad
       e_{\alpha2} = e_1{a}^{32}_{2r}{a}^{42}_{2s}{a}^{34}_{2t-1}, \\
     & e_\beta = \gamma_{12}(e_{\alpha2}) =
        e_1{a}^{32}_{2r}{a}^{42}_{2s}{a}^{34}_{2t}, \text{ where } \\
     & \beta = 1(41)^r(31)^s(21)^t, \quad (\text{case } G11).
   \end{split}
  \end{equation}
    Thus, by Proposition \ref{more_relat_adm}, heading 4), we have
    $\beta = (12)^t(42)^r(32)^s1 = \alpha21$.\vspace{3mm}

 \underline{Case $F32$}. We have
  \begin{equation}
    \begin{split}
     & \alpha2 = (32)^s(42)^r(12)^t, \quad
       e_{\alpha2} = e_3{a}^{12}_{2t}{a}^{42}_{2r}{a}^{14}_{2s-1}, \\
     & e_\beta = \gamma_{12}(e_{\alpha2}) =
        e_3{a}^{12}_{2t+1}{a}^{42}_{2r}{a}^{14}_{2s-1}, \text{ where } \\
     & \beta = 3(41)^s(21)^t(31)^r, \quad (\text{case } G31
       \text{ with } r \leftrightarrow s).
   \end{split}
  \end{equation}
    By Proposition \ref{more_relat_adm}, heading 6), we have
    $\beta = (32)^s(42)^r(12)^t1 = \alpha21$.\vspace{3mm}

  \underline{Case $G22$}. Here,
  \begin{equation}
    \begin{split}
     & \alpha2 = 2(42)^r(32)^s(12)^t, \quad
       e_{\alpha2} = e_2{a}^{14}_{2s}{a}^{34}_{2t}{a}^{31}_{2r}, \\
     & e_\beta = \gamma_{12}(e_{\alpha2}) =
        e_2{a}^{14}_{2s}{a}^{34}_{2t+1}{a}^{31}_{2r}, \text{ where } \\
     & \beta = (21)^{t+1}(41)^s(31)^r, \quad (\text{case } F21
        \text{ with } r \leftrightarrow s, t \rightarrow t+1 ).
   \end{split}
  \end{equation}
    By Proposition \ref{more_relat_adm}, heading 7), we have
    $\beta = (21)^{t+1}(41)^s(31)^r1 = \alpha21$.\vspace{3mm}

 \underline{Case $G32$}. We have
  \begin{equation}
    \begin{split}
     & \alpha2 = 3(42)^r(12)^t(32)^s, \quad
       e_{\alpha2} = e_3{a}^{14}_{2s}{a}^{12}_{2t+1}{a}^{24}_{2r-1}, \\
     & e_\beta = \gamma_{12}(e_{\alpha2}) =
        e_3{a}^{14}_{2s}{a}^{12}_{2t+2}{a}^{24}_{2r-1}, \text{ where } \\
     & \beta = (31)^{r}(41)^s(21)^{t+1}, \quad (\text{case } F31
       \text{ with } r \leftrightarrow s, t \rightarrow t+1).
   \end{split}
  \end{equation}
    By Proposition \ref{more_relat_adm}, heading 8), we have
    $\beta = 3(42)^r(12)^t(32)^s1 =  \alpha21$. \vspace{3mm}

 \underline{Case $G12$}. We have
  \begin{equation}
    \begin{split}
     & \alpha2 = 1(42)^r(32)^s(12)^t, \quad
       e_{\alpha2} = e_1{a}^{34}_{2t}{a}^{32}_{2s+1}{a}^{24}_{2r-1}, \\
     & e_\beta = \gamma_{12}(e_{\alpha2}) =
        e_1{a}^{34}_{2t+1}{a}^{32}_{2s+1}{a}^{24}_{2r-1}, \text{ where } \\
     & \beta = (13)^{r}(41)^{s+1}(21)^t, \quad (\text{case } H11
       \text{ with } r \leftrightarrow s, s \rightarrow s+1).
   \end{split}
  \end{equation}
    By Proposition \ref{more_relat_adm}, heading 9), we have
    $\beta = 1(42)^r(32)^s(12)^t1 =  \alpha21$. \vspace{3mm}

 \underline{Case $H22$}. In this case we have
  \begin{equation}
    \begin{split}
     & \alpha2 = (23)^s(42)^r(12)^t, \quad
       e_{\alpha2} = e_2{a}^{13}_{2r-1}{a}^{14}_{2s-1}{a}^{34}_{2t+1}, \\
     & e_\beta = \gamma_{12}(e_{\alpha2}) =
       e_2{a}^{13}_{2r-1}{a}^{14}_{2s-1}{a}^{34}_{2t+2}, \text{ where } \\
     & \beta = 2(41)^s(31)^{r-1}(21)^{t+1}, \quad (\text{case } G21,
        t \rightarrow t+1, r \leftrightarrow s, r \rightarrow r-1).
   \end{split}
  \end{equation}
    By Proposition \ref{more_relat_adm}, heading 10), we have
    $\beta = (23)^s(42)^r(12)^t1 =  \alpha21$. \vspace{3mm}

   We drop case $F42$ (resp. $G42$) which is similar to
   $F32$ (resp. $G32$) and is proved
   just by permutation $3 \leftrightarrow 4$. \qedsymbol

 \begin{corollary} For cumulative elements $e(n)$, the following relation holds:
   \label{corol_cumulative}
   \begin{equation}
      e_i(n+1) = \gamma_{ij}(e_j(n)) + \gamma_{ik}(e_k(n)) + \gamma_{il}(e_l(n)),
   \end{equation}
 where $\{i,j,k,l\} = \{1,2,3,4\}$.
 \end{corollary}

  \PerfProof
  Without lost of generality, we will show that
  \begin{equation}
    \label{from_en_to_en1}
        e_1(n+1) = \gamma_{12}(e_2(n)) + \gamma_{13}(e_3(n)) + \gamma_{14}(e_4(n)).
  \end{equation}

  By Theorem \ref{relat_rho_ik} we have
\begin{equation*}
  \begin{split}
         & \gamma_{12}(e_2(n)) =
             \sum\limits_{|\alpha| = n-1}e_{\alpha21}, \\
         & \\
         & \gamma_{13}(e_3(n)) =
             \sum\limits_{|\alpha| = n-1}e_{\alpha31}, \\
         & \\
         &  \gamma_{14}(e_4(n)) =
             \sum\limits_{|\alpha| = n-1}e_{\alpha41},
    \end{split}
\end{equation*}
  where $\alpha$ in every sum runs over all admissible elements of
  the length $n-1$. Then,
 \begin{equation*}
      \gamma_{12}(e_2(n)) + \gamma_{13}(e_3(n)) + \gamma_{14}(e_4(n)) =
             \sum\limits_{|\beta| = n}e_{\beta1},
  \end{equation*}
 where $\beta$ runs over all admissible elements of
  the length $n$. The last sum is $e_1(n+1)$ and relation (\ref{from_en_to_en1})
  is proved. \qedsymbol

  Similarly to relation (\ref{lab_rho_ik}) for $e_{\alpha}$,
  it can be proved the relation for admissible elements
  $f_{\alpha0}$:
  \begin{equation}
    \label{prove_rho_2_f}
       \gamma_{ik}(f_{\alpha{k0}}) = f_{\alpha{ki0}}.
  \end{equation}
  For example,
  $$
     f_{20} = e_2(e_1 + e_3 + e_4),
  $$
 and
 \begin{equation*}
   \begin{split}
     \gamma_{12}(f_{20}) = &
        e_2(e_3 + e_4)(e_1(e_3 + e_4) + e_3(e_1 + e_2) + e_4(e_1 + e_2))= \\
     & e_2(e_3 + e_4)(e_1 + e_3(e_1 + e_2) + e_4) = f_{210},
   \end{split}
 \end{equation*}
 see Proposition \ref{coincidence_F} and Section \ref{examples_D4}.

 \subsubsection{A sequence of Herrmann's endomorphisms}
  \label{sect_seq_Herrmann}
   Now, we will see how admissible elements $e_\alpha$ is obtained
   by means of sequence of Herrmann's endomorphisms.
  \begin{theorem}[Admissible elements and Herrmann's endomorphisms]
    \label{theor_adm_Herrmann}
     Admissible elements $e_\alpha$ and $f_{\alpha0}$ ending at $1$
      (from Table \ref{tab_gen_form_D4_e})
     are obtained by means of Herrmann's endomorphisms as follows:
     \begin{equation}
       \label{adm_Herrmann_e}\fbox{$
           \gamma^t_{12}\gamma^s_{13}\gamma^r_{14}(e_1) =
             e_1{a}^{34}_t{a}^{24}_s{a}^{32}_r,$}
     \end{equation}
 and
     \begin{equation}
       \label{adm_Herrmann_f}\fbox{$
           \gamma^t_{12}\gamma^s_{13}\gamma^r_{14}(f_{10}) =
             e_1{a}^{34}_t{a}^{24}_s{a}^{32}_r
              (e_1{a}^{34}_{t+1} + {a}^{24}_{s+1}{a}^{32}_{r-1}),
              \quad \text{ for } r \geq 1.$}
     \end{equation}

     Similarly, admissible elements $e_\alpha$ and $f_{\alpha0}$
     ending at $i = 2,3,4$ are obtained as
     follows:
     \begin{equation}
       \begin{split}
         & \gamma^t_{12}\gamma^s_{13}\gamma^r_{14}(e_2) =
            e_2{a}^{34}_t{a}^{14}_s{a}^{13}_r, \\
         & \gamma^t_{12}\gamma^s_{13}\gamma^r_{14}(e_3) =
            e_3{a}^{14}_t{a}^{24}_s{a}^{12}_r, \\
         & \gamma^t_{12}\gamma^s_{13}\gamma^r_{14}(e_4) =
            e_4{a}^{31}_t{a}^{21}_s{a}^{23}_r,
       \end{split}
     \end{equation}
     and
     \begin{equation}
       \begin{split}
         & \gamma^t_{12}\gamma^s_{13}\gamma^r_{14}(f_{20}) =
             e_2{a}^{34}_t{a}^{14}_s{a}^{31}_r
              (e_2{a}^{34}_{t+1} + {a}^{14}_{s+1}{a}^{31}_{r-1}),
              \quad \text{ for } r \geq 1. \\
         & \gamma^t_{12}\gamma^s_{13}\gamma^r_{14}(f_{30}) =
             e_3{a}^{14}_t{a}^{24}_s{a}^{12}_r
              (e_3{a}^{14}_{t+1} + {a}^{24}_{s+1}{a}^{12}_{r-1}),
              \quad \text{ for } r \geq 1. \\
         & \gamma^t_{12}\gamma^s_{13}\gamma^r_{14}(f_{40}) =
             e_4{a}^{31}_t{a}^{21}_s{a}^{32}_r
              (e_1{a}^{31}_{t+1} + {a}^{21}_{s+1}{a}^{32}_{r-1}),
              \quad \text{ for } r \geq 1.
       \end{split}
     \end{equation}
  \end{theorem}

    \PerfProof 1) Consider the action of Herrmann's endomorphisms on
    $e_\alpha$. For $t=1, s=0, r = 0$ relation (\ref{adm_Herrmann_e})
    follows from (\ref{base_rho12}). Assume,
    relation(\ref{adm_Herrmann_e}) is true for $r+s+t=n$.
    By (\ref{prove_rho_1}) and (\ref{incl_atomic_D4}) we have
  \begin{equation*}
   \begin{split}
    \gamma^{t+1}_{12}\gamma^s_{13}\gamma^r_{14}(e_1) = &
      \gamma_{12}(e_1{a}^{34}_t{a}^{24}_s{a}^{32}_r) =
      \gamma_{12}(e_1{a}^{34}_t) \gamma_{12}(e_1{a}^{24}_s)
      \gamma_{12}(e_1{a}^{23}_r) = \\
    & e_1{a}^{34}_{t+1} e_1{a}^{43}_1a^{24}_s e_1{a}^{34}_1a^{23}_r =
      e_1{a}^{34}_{t+1}a^{24}_s{a}^{23}_r.
   \end{split}
  \end{equation*}
 Thus, (\ref{adm_Herrmann_e}) is true for all triples $(r,s,t)$.

   2) Now, consider action on $f_{\alpha0}$.
   For $t=1, s=0, r = 0$, we have
 \begin{equation}
   \label{rho12_f10}
   \begin{split}
     \gamma_{12}(f_{10}) = & \gamma_{12}(e_1(e_2 + e_3 + e_4)) = \\
    & e_1(e_3 + e_4)
      (e_2(e_3 + e_4) + e_3(e_1 + e_2) + e_4(e_1 + e_2)) = \\
    & e_1{a}^{34}_1(e_2 + e_3(e_1 + e_2) + e_4) =
      e_1{a}^{34}_1(e_2 + e_1(e_3 + e_2) + e_4) = \\
    & e_1{a}^{34}_1(e_1{a}^{32}_1 + a^{42}_1).
   \end{split}
  \end{equation}
   On the other hand, for the case $F12$, $f_{\alpha0}$  can be
   written as follows:
 \begin{equation}
   \label{falpa_t1}
   \begin{split}
      f_{\alpha0} = & e_\alpha(e_1{a}^{34}_{2t} + a^{42}_{2s+1}{a}^{32}_{2r-1})= \\
         & e_\alpha(e_1{a}^{34}_{2t}{a}^{32}_{2r-1} + a^{42}_{2s+1})= \\
         & e_\alpha(e_1{a}^{34}_{2t-2}{a}^{32}_{2r+1} + a^{42}_{2s+1}).
   \end{split}
  \end{equation}
  We use the last expression of $f_{\alpha0}$ in (\ref{falpa_t1}) for the case
  $t = 1, r = 0, s = 0$.
  Then,
 \begin{equation}
   \label{falpa_t1_2}
      f_{\alpha0} =
         e_\alpha(e_1{a}^{34}_{0}{a}^{32}_{1} + a^{42}_{1})=
         e_\alpha(e_1{a}^{32}_{1} + a^{42}_{1}).
  \end{equation}
   Thus, by (\ref{rho12_f10}) and (\ref{falpa_t1_2})
   relation (\ref{adm_Herrmann_f}) is true
   for $t=1, s=0, r = 0$.

   By (\ref{adm_Herrmann_f}) and Proposition \ref{relat_rho_1j}
   we have the induction step:
 \begin{equation*}
   \begin{split}
      \gamma^{t+1}_{12}\gamma^s_{13}\gamma^r_{14}(f_{\alpha0}) = &
        \gamma_{12}(e_1{a}^{34}_t{a}^{24}_s{a}^{32}_r
          (e_1{a}^{34}_{t+1} + {a}^{24}_{s+1}{a}^{32}_{r-1})) = \\
       & \gamma_{12}(e_1{a}^{34}_t)\gamma_{12}(e_1{a}^{24}_s)
         \gamma_{12}(e_1{a}^{32}_r)
         (\gamma_{12}(e_1{a}^{34}_{t+1})
          + \gamma_{12}(e_1{a}^{24}_{s+1})
            \gamma_{12}({a}^{32}_{r-1})) = \\
       & (e_1{a}^{34}_{t+1})(e_1{a}^{34}_1{a}^{24}_s)
         (e_1{a}^{43}_1{a}^{23}_r)
         (e_1{a}^{34}_{t+2} + (e_1{a}^{34}_1{a}^{24}_{s+1})
              (e_1{a}^{43}_1{a}^{23}_{r-1})) = \\
       & e_1{a}^{34}_{t+1}{a}^{24}_s{a}^{23}_r
         (e_1{a}^{34}_{t+2} + e_1{a}^{34}_1{a}^{24}_{s+1}
              e_1{a}^{23}_{r-1}) = \\
       & e_1{a}^{34}_{t+1}{a}^{24}_s{a}^{23}_r
         (e_1{a}^{34}_{t+2} + {a}^{24}_{s+1}e_1{a}^{23}_{r-1}).
         \qed
   \end{split}
  \end{equation*}

 \subsubsection{The sum of Herrmann's endomorphisms}

  Consider endomorphism $\mathcal{R}$ is the sum of Herrmann's endomorphisms
  $\gamma_{1i}$:
  \begin{equation}
    \label{the_sum_R}
     \mathcal{R} = \gamma_{12} + \gamma_{13} + \gamma_{14}.
  \end{equation}

 \begin{proposition}
  \label{prop_primary}
    The endomorphism $\mathcal{R}$ relates
    the inverse cumulative elements (\ref{cumul_polyn_e}) as follows:
  \begin{equation}
    \label{primary_endom}
     \mathcal{R}e_i^\vee(n) = e_i^\vee(n+1), \quad
       \text{ where } i=1,2,3,4.
  \end{equation}
  and
   \begin{equation}
    \label{primary_endom_2}
     \mathcal{R}f_0^\vee(n) = f_0^\vee(n+1).
  \end{equation}

  \end{proposition}

  \PerfProof
   1) By (\ref{cumul_polyn_e})
 \begin{equation*}
    e_1^\vee(n) = \sum\limits_{r + s + t = n-1}
                      {e_1{a}_r^{32}{a}_s^{24}{a}_t^{34}},
 \end{equation*}
 and by Proposition \ref{relat_rho_1j} we have
 \begin{equation*}
   \begin{split}
     \mathcal{R}(e_1^\vee(n))
         = & (\gamma_{12} + \gamma_{13} + \gamma_{14})\sum\limits_{r + s + t = n-1}
                      {e_1{a}_r^{32}{a}_s^{24}{a}_t^{34}} = \\
           & \sum\limits_{r + s + t = n-1}{e_1{a}_r^{32}{a}_s^{24}{a}_{t+1}^{34}} +
             \sum\limits_{r + s + t = n-1}{e_1{a}_r^{32}{a}_{s+1}^{24}{a}_t^{34}} +
             \sum\limits_{r + s + t = n-1}{e_1{a}_{r+1}^{32}{a}_s^{24}{a}_t^{34}}
             = \\
           & \sum\limits_{r + s + t = n}
                      {e_1{a}_r^{32}{a}_s^{24}{a}_t^{34}} =
                      e_1^\vee(n+1).
   \end{split}
 \end{equation*}

   2) By (\ref{cumul_polyn_f})
 \begin{equation*}
    f_0^\vee(n) = \sum\limits_{r + s + t = n}
                      {e_\alpha({a}_r^{jk} + {a}_s^{kl}{a}_t^{lj}}),
 \end{equation*}
 and, again, by Proposition \ref{relat_rho_1j}
 \begin{equation*}
   \begin{split}
     \mathcal{R}(f_0^\vee(n)) =
           & (\gamma_{12} + \gamma_{13} + \gamma_{14})\sum\limits_{r + s + t = n}
                  e_\alpha({a}_r^{jk} + {a}_s^{kl}{a}_t^{lj}) = \\
           & (\gamma_{12} + \gamma_{13} + \gamma_{14})\sum\limits_{r + s + t = n}
                  \gamma^t_{12}\gamma^s_{23}\gamma^r_{14}(f_{i0}), \\
           &      \text {where } \{i,j,k,l\} = \{1,2,3,4\}.
   \end{split}
 \end{equation*}
 Thus,
 \begin{equation*}
   \begin{split}
    \mathcal{R}(f_0^\vee(n)) =
             & \sum\limits_{r + s + t = n}
                \gamma^{t+1}_{12}\gamma^s_{23}\gamma^r_{14}(f_{i0}) + \\
             & \sum\limits_{r + s + t = n}
                \gamma^{t}_{12}\gamma^{s+1}_{23}\gamma^r_{14}(f_{i0}) +
             \sum\limits_{r + s + t = n}
                \gamma^{t}_{12}\gamma^{s}_{23}\gamma^{r+1}_{14}(f_{i0}) = \\
       & \sum\limits_{r + s + t = n+1}
                \gamma^{t}_{12}\gamma^s_{23}\gamma^r_{14}(f_{i0}) =
                f_0^\vee(n+1). \qed
  \end{split}
 \end{equation*}

  \subsection{Perfect elements in $D^4$}
  \subsubsection{The Gelfand-Ponomarev conjecture}
    \label{sect_GP_conjecture}

 By Gelfand-Ponomarev definition (\ref{perfect_D4}), elements
 $h_i(n)$, where $i = 1,2,3,4$  are generators of the $16$-element
 Boolean cube $B^+(n)$ of perfect elements.
  According to (\ref{perfect_D4}), these $16$ elements in $B^+(n)$
  are as follows:
 \begin{equation}
   \label{elem_h_D4_max_min}
   \begin{split}
      & h^{\max}(n) = \sum\limits_{i=1,2,3,4}h_i(n) = \sum\limits_{i=1,2,3,4}e_i(n), \\
      & \\
      & h^{\min}(n) = \bigcap\limits_{i = 1,2,3,4}h_i(n) =
          \sum\limits_{i=1,2,3,4}e_i(n)h_i(n),
   \end{split}
 \end{equation}
 and
 \begin{equation}
   \label{elem_h_D4}
   \begin{split}
      & h_i(n) = e(j) + e(k) + e(l), \\
      & h_i(n)h_j(n) = e_i(n)h_i(n) + e_j(n)h_j(n) + e_k(n) + e_l(n), \\
      & h_i(n)h_j(n)h_k(n) = e_i(n)h_i(n) + e_j(n)h_j(n) + e_k(n)h_k(n) + e_l(n), \\
      &  \qquad \qquad \text{ where } \{i,j,k,l\} = \{1,2,3,4\},
   \end{split}
 \end{equation}
 see also Table \ref{tab_perfect_in D4}.

 Four elements $h_i(n)$ (resp. four elements $h_i(n)h_j(n)h_k(n)$)
 from (\ref{elem_h_D4}) are coatoms (resp. atoms) in $B^+(n)$. The
 maximal and minimal elements in $B^+(1)$ are as follows:
 \begin{equation}
   \label{perfect_D4_max_min_ex}
   \begin{split}
     & h^{\max}(1) = e_1 + e_2 + e_3 + e_4, \\
     & h^{\min}(1) =
       (e_1 + e_2 + e_3)(e_1 + e_2 + e_4)(e_1 + e_3 + e_4)(e_2 + e_3 + e_4). \\
   \end{split}
 \end{equation}

  By (\ref{perfect_D4_max_min_ex}) and Section \ref{examples_D4} we have
 \begin{equation}
   \label{perfect_D4_min_1}
   \begin{split}
       h^{\min}(1) =
        & e_1(e_2 + e_3 + e_4) + e_2(e_1 + e_3 + e_4) + e_3(e_1 + e_2 + e_4) =\\
        & f_{10} + f_{20} + f_{30}.
   \end{split}
 \end{equation}
 By (\ref{perfect_D4_max_min_ex}) also
 $$
    h^{\min}(1) \supset e_4(e_1 + e_2 + e_3) = f_{40},
 $$
 and by (\ref{perfect_D4_min_1}) we obtain
 $$
    h^{\min}(1) = f_{10} + f_{20} + f_{30} + f_{40},
 $$
 Thus,
 \begin{equation}
   \label{min_elem_in_B1}
    h^{\min}(1) = f_0(2).
 \end{equation}
 Relation (\ref{min_elem_in_B1}) takes place modulo linear
 equivalence for any $n$:
 \begin{equation}
   \label{min_elem_in_Bn}
    h^{\min}(n) \simeq f_0(n+1).
 \end{equation}
 Relation (\ref{min_elem_in_Bn}) is proved like the similar
 relation for $D^{2,2,2}$, see \cite[Prop. 3.4.2]{St04}.
  According to the Gelfand-Ponomarev conjecture
 \cite[p.7]{GP74}, the relation (\ref{min_elem_in_Bn})
 takes place:
 \begin{equation}
   \label{GP_conjecture}
    h^{\min}(n) = f_0(n+1) \quad
      \text{(Gelfand-Ponomarev conjecture)}.
 \end{equation}
 Further, by (\ref{elem_h_D4})
 \begin{equation}
   \label{h1_h2_h3_1}
   \begin{split}
    h_1(1)h_2(1)h_3(1) = & (e_2 + e_3 + e_4)(e_1 + e_3 + e_4)(e_1 + e_2 + e_4)= \\
          & e_1(e_2 + e_3 + e_4) + e_2(e_1 + e_3 + e_4) + e_3(e_1 + e_2 +  e_4) + e_4 = \\
          & f_{10} + f_{20} + f_{30} +  e_4 =  e_4 + f_0(2).
   \end{split}
 \end{equation}

 \begin{proposition}
     Elements $h_i(n)h_j(n)h_k(n)$ and $e_l(n) + f_0(n+1)$
    coincide modulo linear equivalence for all $\{i,j,k,l\} =
    \{1,2,3,4\}$:
   \begin{equation}
     \label{h_ijk_n}
      h_i(n)h_j(n)h_k(n) \simeq e_l(n) + f_0(n+1).
   \end{equation}
 \end{proposition}
  \PerfProof
  Assuming that (\ref{h_ijk_n}) is true, we will prove it for  $n+1$:
 \begin{equation}
   \label{h1_h2_h3_n1}
   \begin{split}
    h_i(n+1)h_j(n+1)h_k(n+1) \simeq e_l(n+1) + f_0(n+2).
   \end{split}
 \end{equation}
  Similarly to \cite[Prop. 3.4.1]{St04} for $D^{2,2,2}$,
  we have the following relation in $D^4$ for
  sums and intersections, commuting on the perfect elements $v_t$:
 \begin{equation}
   \label{sum_inters_D4}
    \sum_\text{$p = 1,2,3,4$}
      \varphi_p\rho_{X^+}(\bigcap_\text{$t=i,j,k$}v_t) =
        \bigcap_\text{$t=i,j,k$}(\sum_\text{$p =
        1,2,3,4$}\varphi_p\rho_{X^+}(v_t)),
 \end{equation}
 where all $v_t$ are perfect elements\footnote{
  In the proof of Proposition $3.4.1.$ from \cite{St04} we only
  change $\sum\varphi_i\rho_{X^+}(I)$.
  According to Corollary \ref{cor_psi_D4}, heading (3), in the case of $D^4$ we
  have
 \begin{equation*}
    \sum\limits_{i=1,2,3,4}\varphi_i\rho_{X^+}(I) =
      \sum\limits_{i=1,2,3,4}\rho_X(e_i(e_j + e_k + e_l)) =
      \rho_X(\bigcap\limits_{i=1,2,3,4}(e_j + e_k + e_l)).
 \end{equation*}
 }. 
Since
 $$
    \sum_\text{$p = 1,2,3,4$}\varphi_p\rho_{X^+}(h_t(n)) = \rho_X(h_t(n+1)),
 $$
   we have in the right side of (\ref{sum_inters_D4}), that
 \begin{equation}
  \label{r_side}
   \bigcap_\text{$t=i,j,k$}(\sum_\text{$p = 1,2,3,4$}\varphi_p\rho_{X^+}(v_t)) =
     \rho_X(\bigcap\limits_{t=i,j,k}h_t(n+1)).
 \end{equation}
 In the left side of (\ref{sum_inters_D4}) we have
 \begin{equation}
  \label{l_side}
  \begin{split}
   \sum_\text{$p = 1,2,3,4$}
      \varphi_p\rho_{X^+}(\bigcap_\text{$t=i,j,k$}h_t(n)) = &
   \sum_\text{$p = 1,2,3,4$}
      \varphi_p\rho_{X^+}(e_l(n) + f_0(n+1)) = \\
      & \\
      & \rho_X(e_l(n+1) + f_0(n+2)).
   \end{split}
 \end{equation}
 Relation (\ref{h1_h2_h3_n1}) follows from (\ref{r_side}) and (\ref{l_side}).
 \qedsymbol

 Similarly, other relations from Table \ref{tab_perfect_in D4} take
 place modulo linear equivalence:
 \begin{equation}
   \label{hi_hj_hk_n}
   \begin{split}
    & h_i(n)h_j(n)h_k(n) \simeq e_l(n) + f_0(n+1), \\
    & h_i(n)h_j(n) \simeq e_k(n) +  e_l(n) + f_0(n+1), \\
    & h_i(n) \simeq e_j(n) + e_k(n) +  e_l(n) + f_0(n+1).
   \end{split}
 \end{equation}

 \begin{figure}[h]
 \includegraphics{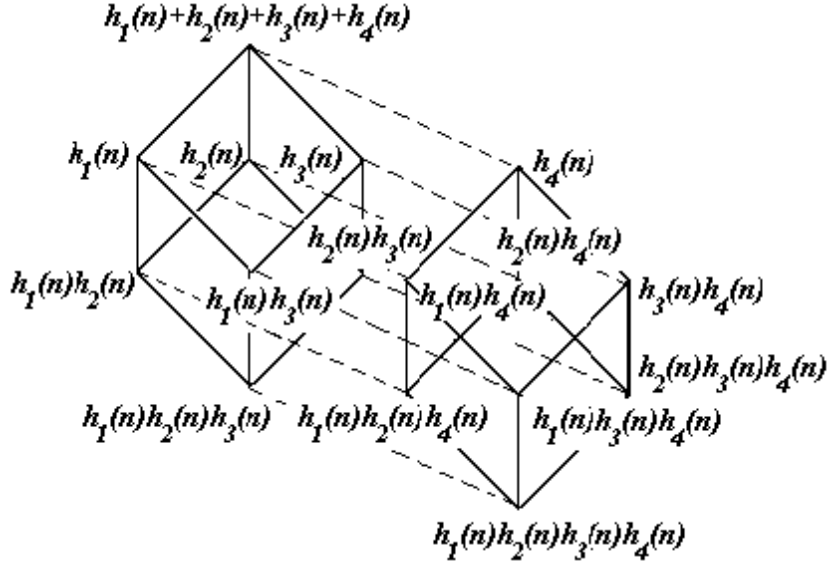}
 \caption{\hspace{3mm}The $16$-element Boolean cube $B^+(n)$ with
 generators $h_i(n)$}
 \label{h_Bn_cube}
 \end{figure}

 \begin{proposition}
    The following relations hold:
 \begin{equation}
    \begin{split}
      \label{sum_h_ijk_ijl}
      &  \quad h_i(n)h_j(n)h_k(n) + h_i(n)h_j(n)h_l(n) = h_i(n)h_j(n),\\
      &  \quad h_i(n)h_j(n) + h_i(n)h_k(n) = h_i(n).
    \end{split}
  \end{equation}
 \end{proposition}
 \PerfProof
   It easily follows from (\ref{elem_h_D4}) and
   (\ref{perfect_D4}).
   \qedsymbol

  \subsubsection{Herrmann's polynomials $s_n$, $t_n$ and $p_{i,n}$}
   \label{sect_stp_poly}
  C.~Herrmann used endomorphisms $\gamma_{ij}$ from (\ref{endom_rho_ij})
  for the construction of perfect elements $s_n$, $t_n$ and
  $p_{i,n}$, where $i = 1,2,3,4$, see \cite[p. 362]{H82}, \cite[p. 229]{H84}.
  In what follows, the definitions of

  \underline{Polynomials $s_n$:}
  \begin{equation}
   \label{def_sn}
   \begin{split}
      & s_0 = I, \quad s_1 = e_1 + e_2 + e_3 + e_4, \\
      & s^{i-1}_n = \gamma_{1i}(s_n), \text{ where } i = 2,3,4, \\
      & s_{n+1} = s^1_n + s^2_n + s^3_n.
   \end{split}
  \end{equation}

  \underline{Polynomials $t_n$:}
  \begin{equation}
   \label{def_tn}
   \begin{split}
      & t_0 = I, \quad
        t_1 = (e_1 + e_2 + e_3)(e_1 + e_2 + e_4)(e_1 + e_3 + e_4)(e_2 + e_3 + e_4), \\
      & t^{i-1}_n = \gamma_{1i}(t_n), \text{ where } i = 2,3,4, \\
      & t_{n+1} = t^1_n + t^2_n + t^3_n.
   \end{split}
  \end{equation}

  \underline{Polynomials $p_{k,n}$:}
  \begin{equation}
   \label{def_pnk}
   \begin{split}
      & p_{i,0} = I, \quad
        p_{i,1} = e_i + t_1, \text{ where } i = 1,2,3,4, \\
      & p_{i,n+1} = \gamma_{ij}(p_{j,n}) + \gamma_{ik}(p_{k,n}) + \gamma_{il}(p_{l,n}), \\
      &\qquad \qquad
        \text{ where } i = 1,2,3,4, \text{ and } \{i,j,k,l\} = \{1,2,3,4\}.
   \end{split}
  \end{equation}
  For example, by (\ref{the_sum_R}) we have
  \begin{equation}
   \label{def_pnk_example}
   \begin{split}
      p_{1,2} = & \gamma_{12}(p_{2,1}) + \gamma_{13}(p_{3,1}) + \gamma_{14}(p_{4,1}) =
         \gamma_{12}(e_2) + \gamma_{13}(e_3) +  \gamma_{14}(e_4) + \mathcal{R}(t_1)= \\
      & e_2(e_3 + e_4) + e_3(e_2 + e_4) + e_4(e_2 + e_3) + t_2 =
        e_{21} + e_{31} + e_{41} + t_2 = \\
      & e_1(2) + t_2.
   \end{split}
  \end{equation}
 Similarly,
 $$
     p_{i,2} = e_i(2) + t_2.
 $$

\begin{figure}[h]
 \includegraphics{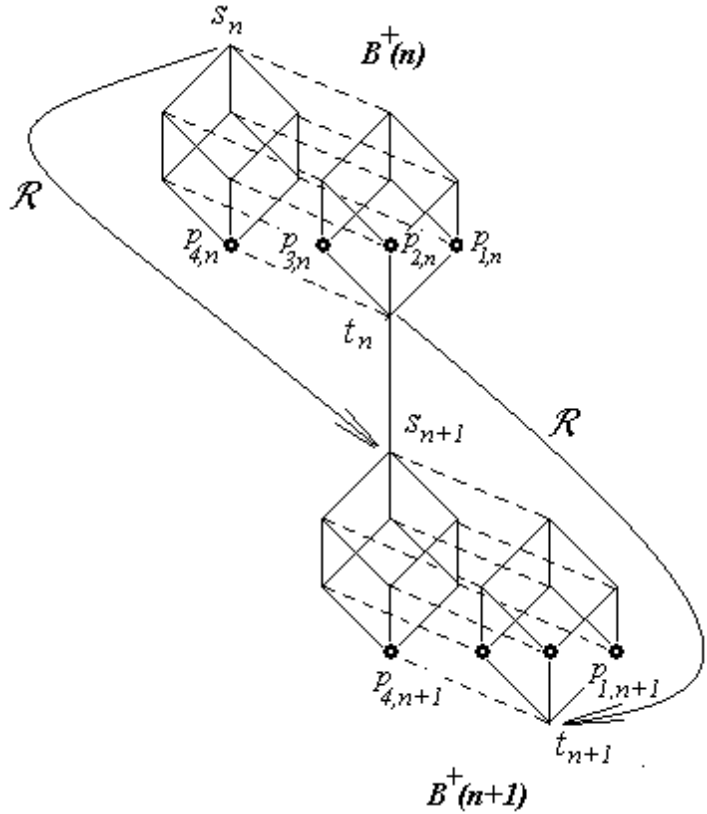}
 \caption{\hspace{3mm}The endomorphism $\mathcal{R}$ mapping $s_n$
 to $s_{n+1}$ and $t_n$ to $t_{n+1}$}
 \label{endom_R}
 \end{figure}

  By definitions (\ref{def_sn}), (\ref{def_tn}) we have
  \begin{equation}
    \begin{split}
      & s_{n+1} = \mathcal{R}s_n = \mathcal{R}^{n}s_1, \\
      & t_{n+1} = \mathcal{R}t_n = \mathcal{R}^{n}t_1, \\
     \end{split}
  \end{equation}
 see Fig. \ref{endom_R}.

 \subsubsection{Cumulative elements and Herrmann's polynomials}
 In what follows, we show how Herrmann's polynomials $s_n$,
 $t_n$ and $p_{i,n}$, where $i \in \{1,2,3,4\}$, are calculated by means
 of cumulative elements $e(n)$, $f(n)$.
 Polynomials $s_n$ and $t_n$ can be also calculated by means of
 inverse cumulative elements $e^\vee(n)$, $f^\vee(n)$,
 see Section \ref{sect_inverse_cum}. As above, the main elementary brick in these
 constructions is an admissible element, obtained as a
 sequence of Herrmann's endomorphisms
 $$
    \gamma^r_{12}\gamma^s_{13}\gamma^t_{14}(e_i),
 $$
 see Theorem \ref{theor_adm_Herrmann} from Section \ref{sect_seq_Herrmann}.

 \begin{theorem}[The polynomials $s_n$, $t_n$ and $p_{i,n}$]
  \label{prop_endom_R}
     1) For each $n$, the polynomial $s_n$ is the maximal perfect element in the
      Boolean cube $B^+(n)$, see Section \ref{sect_perfect_D4},
      and it is expressed as follows:
    \begin{equation}
     \label{intro_polyn_Herrmann}
      \begin{split}
       s_n = \mathcal{R}^{n-1}(s_1) = &
        \sum\limits_{r + s + t = n-1}
        \gamma^t_{12}\gamma^s_{13}\gamma^r_{14}(e_1 + e_2 + e_3 + e_4) = \\
        & \\
        & \sum\limits_{i=1,2,3,4}e_i^\vee(n) = \sum\limits_{i=1,2,3,4}e_i(n). \\
        & \\
      \end{split}
    \end{equation}

     2) For each $n$, the polynomial $t_n$ is linearly equivalent to
      the minimal perfect element in the
      Boolean cube $B^+(n)$, see Section \ref{sect_perfect_D4},
      and it is expressed as follows:
    \begin{equation}
     \label{tn_minimal}
      \begin{split}
        & \\
        t_n = & \mathcal{R}^{n-1}(t_1) = \\
              & \\
              & \sum\limits_{r + s + t = n}
         \gamma^t_{12}\gamma^s_{13} \gamma^r_{14}
          (e_1(e_2 + e_3 + e_4) + e_2(e_1 + e_3 + e_4) +
            e_3(e_1 + e_2  + e_4)) = \\
            & \\
            & f^\vee_0(n+1) = f_0(n+1). \\
            & \\
      \end{split}
    \end{equation}

     3) For each $n$, the polynomial $p_{i,n}$ is linearly equivalent to
      the element $h_j(n)h_k(n)h_l(n)$ in the
      Boolean cube $B^+(n)$, see Table \ref{tab_perfect_in D4},
      and it is expressed as follows:
    \begin{equation}
     \label{p_1n_minimal}
        p_{i,n} = e_i(n) + f^\vee_0(n+1) = e_i(n) + f_0(n+1).
    \end{equation}
 \end{theorem}

  \PerfProof 1) By definition (\ref{def_sn}) of $s_n$, and
    by property (\ref{prove_rho_1}),  for $n=1$, we have
 \begin{equation}
  \label{calc_s2}
   \begin{split}
      & s_2 = \gamma_{12}(s_1) + \gamma_{13}(s_1) + \gamma_{14}(s_1) =
        \sum\limits_{i=1,2,3,4}
        (\gamma_{12} + \gamma_{13} + \gamma_{14})(e_i) = \\
      & \sum\limits_{\substack{j=1,2,3 \\ i=1,2,3,4}}\gamma_{1j}(e_i)=
        \sum\limits_{\substack{i=1,2,3,4 \\ k,l \neq i}}e_i{a}^{kl}_1 =
        \sum\limits_{|\alpha| = 2}e_{\alpha} =
        e_1^\vee(2) + e_2^\vee(2) + e_3^\vee(2) + e_4^\vee(2).
   \end{split}
  \end{equation}
  In (\ref{calc_s2}) $|\alpha|$ is the length of sequence $\alpha$.
  Note, that in (\ref{calc_s2}) $|\alpha| = 2$ and $r+s+t =1$.
  Assume,
 \begin{equation}
      s_n = \sum\limits_{|\alpha| = n}e_{\alpha} =
        e_1^\vee(n) + e_2^\vee(n) + e_3^\vee(n) + e_4^\vee(n) ,
  \end{equation}
 Then, by definition (\ref{def_sn}) by
   Proposition \ref{prop_primary} and Proposition \ref{sum_of_all} we have
 \begin{equation}
   \begin{split}
     s_{n+1} =
       & (\gamma_{12} + \gamma_{13} + \gamma_{14})\sum\limits_{|\alpha| = n}e_{\alpha} =
         (\gamma_{12} + \gamma_{13} + \gamma_{14})
        \sum\limits_{\substack{i=1,2,3,4 \\ r + s+ t =n-1}}
        \gamma^r_{12}\gamma^s_{13}\gamma^t_{14}(e_i) = \\
       & \\
       &  \sum\limits_{\substack{i=1,2,3,4 \\ r + s + t =n}}
        \gamma^r_{12}\gamma^s_{13}\gamma^t_{14}(e_i) =
         \sum\limits_{|\alpha| = n+1}e_{\alpha} =
         \sum\limits_{i=1,2,3,4}e_i^\vee(n+1) =
         \sum\limits_{i=1,2,3,4}e_i(n+1).
   \end{split}
  \end{equation}
  and $s_{n+1}$ is the maximal element in $B^+(n+1)$, see
  Section \ref{sect_perfect_D4}. \vspace{3mm}

  2) It is easy to see that $t_1$ in (\ref{def_tn}) is as follows:
 \begin{equation}
  \begin{split}
    t_1 = & e_1(e_2 + e_3 + e_4) + e_2(e_1 + e_3 + e_4) +
            e_3(e_1 + e_2  + e_4) = \\
          & f_{10} + f_{20} + f_{30} = f_0(2),
   \end{split}
  \end{equation}
  and by (\ref{adm_Herrmann_f})
 from Theorem \ref{theor_adm_Herrmann} we have
 \begin{equation}
  \label{t2_f03}
  \begin{split}
    t_2 = & \mathcal{R}(f_0(2)) =
      (\gamma_{12} + \gamma_{13} + \gamma_{14})
      \sum\limits_{i=1,2,3,4}f_{i0} =
      \sum\limits_{i \neq j}f_{ij0} = \\
      & \sum\limits_{|\alpha| = 2}f_{\alpha0} = f_0(3).
   \end{split}
  \end{equation}
 Assume, (\ref{tn_minimal}) is true for some $n$, then by
 (\ref{cumul_polyn_e}) we have
 \begin{equation}
 \begin{split}
    t_n = \sum\limits_{\substack{i=1,2,3,4 \\ r + s + t = n}}
        \gamma^r_{12}\gamma^s_{13}\gamma^t_{14}(f_{i0}) = f^\vee_0(n+1).
   \end{split}
  \end{equation}
 Then, by Proposition \ref{prop_primary} and Theorem
 \ref{theor_adm_Herrmann} we have
 \begin{equation}
 \begin{split}
    t_{n+1} = Rf^\vee_0(n) = f^\vee_0(n+1) =
      \sum\limits_{\substack{i=1,2,3,4 \\ r + s + t = n+1}}
        \gamma^r_{12}\gamma^s_{13}\gamma^t_{14}(f_{i0}).
   \end{split}
  \end{equation}
   By Proposition \ref{sum_of_all} $f^\vee_0(n) = f(n)$ and
   according to (\ref{min_elem_in_Bn}), the element $t_{n+1} = f_0(n+1)$ is linearly
   equivalent to the minimal element in $B^+(n+1)$.
   If the Gelfand-Ponomarev conjecture (\ref{GP_conjecture}) is true,
   the element $t_{n+1} = f_0(n+1)$ coincides with
   the minimal element $h^{\min}(n+1)$ of $B^+(n+1)$, see Section \ref{sect_GP_conjecture}.

   3) For n = 1,  $p_{i,1} = e_i + t_1 = e_i(1) + f_0(2)$, and
 \begin{equation}
  \begin{split}
     p_{i,2} =
      & \gamma_{ji}(e_j + t_1) + \gamma_{ki}(e_k + t_1) + \gamma_{li}(e_l + t_1) = \\
      & e_i(2) + \mathcal{R}t_1 = e_i(2) + t_2 = e_i(2) + f_0(3),
   \end{split}
 \end{equation}
   see (\ref{def_pnk_example}), (\ref{t2_f03}).
   Further, by Corollary \ref{corol_cumulative} we have
 \begin{equation}
 \begin{split}
    p_{i,n+1} = & \gamma_{ji}(p_{j,n}) + \gamma_{ki}(p_{k,n}) + \gamma_{li}(p_{l,n}) = \\
    & \gamma_{ji}(e_j(n) + t_n) + \gamma_{ki}(e_k(n) + t_n) + \gamma_{li}(e_l(n) + t_n) = \\
    &  \gamma_{ji}(e_j(n)) + \gamma_{ki}(e_k(n)) + \gamma_{li}(e_l(n)) +  R(t_n) =
              e_i(n+1) + t_{n+1}.
   \end{split}
  \end{equation}
    By heading 2)
 \begin{equation}
    p_{i,n+1} = e_i(n+1) + t_{n+1} = e_i(n+1) + f_0(n+1).
  \end{equation}
    The theorem is proved. \qedsymbol

As corollary we obtain the following
\begin{theorem}[C.~Herrmann, \cite{H82}]
  \label{th_Herrmann_82}
    Polynomials $p_{1,n}, p_{2,n}, p_{3,n}, p_{4,n}$ are atoms in
    $16$-element Boolean cube $B^+(n)$ of perfect elements and
 \begin{equation}
  \label{lab_pn}
     p_{i,n} \simeq  h_j(n)h_k(n)h_l(n), \quad \text{ where } \{i,j,k,l\} = \{1,2,3,4\}, \\
  \end{equation}
     Polynomials $s_n$ and $t_n$ are, respectively, the maximal
     and minimal elements in $B^+(n)$, and
 \begin{equation}
   \label{lab sn_tn}
     s_n \simeq  \sum\limits_{i=1,2,3,4}h_j(n), \quad
     t_n \simeq \bigcap\limits_{1=1,2,3,4}h_j(n),\\
  \end{equation}
see Table \ref{tab_perfect_in D4}.
\end{theorem}
 \PerfProof
  From (\ref{p_1n_minimal}) and (\ref{hi_hj_hk_n})
  we obtain (\ref{lab_pn}). From (\ref{sn_as_sum}), (\ref{tn_minimal})
  and (\ref{elem_h_D4_max_min})
  we obtain (\ref{lab sn_tn}). \qedsymbol